\documentclass[a4paper,12pt]{amsart}
\usepackage{amsmath,amsfonts,amssymb,multicol, amsthm}
\usepackage{amsrefs}
\usepackage{mathrsfs}
\usepackage[utf8]{inputenc}
\usepackage{esint}
\usepackage{color, framed}

\usepackage{graphicx}
\usepackage[normalem]{ulem}

\usepackage{hyperref}
\usepackage{bm}
\usepackage{mathtools}
\usepackage{graphicx}
\usepackage{amsmath,amssymb,amsthm,amsfonts}
\usepackage{amssymb}
\usepackage[active]{srcltx}
\usepackage{bbold}
\usepackage{color}
\newtheorem{thm}{Theorem}[section]
\newtheorem{cor}[thm]{Corollary}
\newtheorem{lem}[thm]{Lemma}
\newtheorem{prop}[thm]{Proposition}
\newtheorem{rem}[thm]{Remark}

\newtheorem*{theorem*}{Theorem}

\usepackage{tikz}
\usepackage{standalone}
\usepackage{graphicx}
\usepackage{subcaption}
\usepackage{multicol,float}

\newcommand{\disp}{\displaystyle}

\numberwithin{equation}{section}

\newcommand{\R}{\mathbb R}

\begin{document}
	
	\parindent 0pc
	\parskip 6pt
	\overfullrule=0pt
	
	\title{\quad One-dimensional symmetry results \\ for semilinear equations and inequalities \\ on half-spaces}
	
	\author{Nicolas Beuvin$^{\dagger}$, Alberto Farina$^{\dagger}$}
	\address{$^{\dagger}$  Universit\'e de Picardie Jules Verne, LAMFA, CNRS UMR 7352, 33, rue Saint-Leu 80039 Amiens, France}\email{nicolas.beuvin@u-picardie.fr, alberto.farina@u-picardie.fr}

	
	

	\maketitle

	\begin{abstract}
We prove new one-dimensional symmetry results for non-negative solutions, possibly unbounded, to the semilinear equation 
$ -\Delta u= f(u)$ in the upper half-space $\R^{N}_{+}$. \\
Some Liouville-type theorems are also proven in the case of differential inequalities in $\R^{N}_{+}$, even without imposing any boundary condition. \\
Although subject to dimensional restrictions, our results apply to a broad family of functions $f$. In particular, they apply to all non-negative $f$ that behaves at least linearly at infinity. 
    \end{abstract}

\section{Introduction}
In this work we consider solutions, not necessarily bounded, to the elliptic boundary value problem
\begin{equation}\label{Poisson}
\left\{
	\begin{array}{ccl}
		-\Delta u= f(u) & \text{in}&  \R^{N}_{+},\\
		u\geq 0 & \text{in} & \R^{N}_{+},\\
		u=0 & \text{on} & \partial \R^{N}_{+},
	\end{array}
	\right.
\end{equation}
where $\mathbb{R}^N_+ =\{x=(x', x_N) \in \R^{N-1} \times \R  \ | \ x_N>0\}$ is the open upper half-space of $\R^N$, $N \geq 2$
and $f$ is a locally Lipschitz-continuous function on $[0,+\infty)$. \\
Studying and classifying solutions to problem \eqref{Poisson} is natural for various reasons:
\begin{enumerate}
\item [(i)] the half-space is the simplest unbounded domain with an unbounded boundary in which one can study the Dirichlet problem 
    for the semilinear Poisson equation. 

\item [(ii)] Problem \eqref{Poisson} arises as a limit problem, after a blow-up procedure near the boundary, either when one wants 
to derive \textit{a priori} estimates for positive solutions to some subcritical semilinear uniformly elliptic BVPs on smooth bounded domains (\cite{GS}), or when one studies uniqueness results for classical solutions to the singular perturbation problem 
\begin{equation}\label{eq1}
- \varepsilon^2 \Delta v = f(v)  \quad \text{in} \quad \Omega, \qquad  v > 0 \quad \text{in} \quad \Omega, \quad v = 0 \quad \text{on} \quad \partial \Omega,
\end{equation}
where $\Omega$ is a smooth and bounded domain, $f$ is a suitable smooth function and $\varepsilon$ is a small and positive parameter (see for instance \cite{An} and \cite{dan2}). 

Notice that, when $f(0) <0$, even if solutions $v$ considered in \eqref{eq1} are positive, given the absence of the strong maximum principle, the solution $u$ of \eqref{Poisson} obtained after the scaling procedure is "only" non-negative in $\R^{N}_{+}$. 
This motivates the choice to consider the non-negative solutions of \eqref{Poisson} (and not only the positive ones).
\end{enumerate}
Finally, the need to also consider unbounded solutions to \eqref{Poisson} is well illustrated by the case of harmonic functions, that is when $f \equiv 0$. In this case, it is well-known that all the solutions to \eqref{Poisson} are given by $u(x)= \alpha x_N$, $ \alpha \geq 0$. 
\smallskip

The first result concerning the complete classification of the solutions to \eqref{Poisson} with $f \not \equiv 0 $ was obtained by Gidas and Spruck in 1981. It concerns the function $f(u) =u^p$ and it was motivated by the the first situation described in item (ii) above.

\textbf{Theorem A} (\cite{GS}) \textit{
Assume $N \geq 2$ and let $p$ be any real number satisfying $ 1 < p \leq \frac{N+2}{N-2}$. The only classical solution to
\begin{equation*}
\left\{
	\begin{array}{ccl}
		-\Delta u= u^p& \text{in}&  \R^{N}_{+},\\
		u \geq 0 & \text{in} & \R^{N}_{+},\\
		u=0 & \text{on} & \partial \R^{N}_{+},
	\end{array}
	\right.
\end{equation*}	
is $ u \equiv 0$. }

The proof in \cite{GS} makes use of the moving spheres method. In 1993, Berestycki, Caffarelli and Nirenberg (\cite{BCNMagenes}) observed that the proof in \cite{GS} actually yields the following result 

\textbf{Theorem B} (\cite{BCNMagenes}) \textit{
Let $u \in C^{2}(\overline{\R^{N}_{+}})$ be a solution of \eqref{Poisson} where 
$f \in Lip_{\text{loc}}([0,+\infty))$. \\
If any of the two following conditions holds:  
\begin{enumerate}
   \item [(i)] $N=2$ and $f \geq 0$,
   \item [(ii)] $ N \geq 3$ and $f$ satisfies 
\begin{equation}\label{cond-crit-gen}
t \mapsto \frac{f(t)}{t^{(N+2)/(N-2)}} \quad \text{is nonincreasing on} \,\, (0,+\infty).\\
\end{equation}
\end{enumerate}
Then, either $u \equiv 0$, or there exists a function $u_{0}:[0,+\infty) \to [0,+\infty)$ such that
	\begin{equation}\label{Def-1D}
		u(x)=u_{0}(x_N) \qquad \text{for any } x \in \R^{N}_{+}
	\end{equation}
and $u_0>0$, $u'_0>0$ on $(0,+\infty)$. 
}

Solutions of the form \eqref{Def-1D}, that is, solutions depending only on the variable $x_N$, are called \textit{one-dimensional} since they inherit the symmetry of the underlying domain $\R^{N}_{+}$. Clearly, the profile $u_0$ solves the one-dimensional problem  \begin{equation}\label{1D-Poisson}
\left\{
	\begin{array}{ccl}
		- {u_0}{''} = f(u_0) & \text{in} & (0,+\infty),\\
		u_0(0)=0. & \text{}
	\end{array}
\right.
\end{equation} 

Note that we necessarily have $f(0) \geq 0$ in the above Theorem B. 
Recently, the following results about the one-dimensional symmetry of the solutions to \eqref{Poisson} have been obtained. They also cover the case in which $f(0)<0$.  

\textbf{Theorem C} (\cite{fs1,fs2}) \textit{
Let $u \in C^{2}(\overline{\R^{2}_{+}})$ be a solution of \eqref{Poisson} where 
$f \in Lip_{\text{loc}}([0,+\infty))$. Then, 
\begin{enumerate}
    \item [(i)] if $f(0)\geq 0$, then either $u \equiv 0$, or $u$ is positive on $\R^{2}_{+}$ and strictly increasing in the direction orthogonal to the boundary, i.e., $\frac{\partial u}{\partial x_2}> 0$ on $\R^{2}_{+}$;
   \item [(ii)] if $f(0) < 0$, then either $u>0$ and $\frac{\partial u}{\partial x_2}> 0$ on $\R^{2}_{+}$, 
   or 
\begin{equation*}
		u(x)=u_{0}(x_2) \qquad \text{for any } x \in \R^{2}_{+}
	\end{equation*}   
for a unique, periodic function $u_{0}:[0,+\infty) \to [0,+\infty)$;
    \item [(iii)] if for some $p>1$, 
    \begin{equation}\label{asym-tp}
	\lim\limits_{t \to + \infty} \frac{f(t)}{t^p} \in (0, +\infty),
\end{equation}
then $u$ is bounded and one-dimensional.\footnote{\, more precisely, $u(x)=u_{0}(x_2)$ for any $ x \in \R^{2}_+$ and for some bounded function $u_{0}:[0,+\infty) \to [0,+\infty)$ solving \eqref{1D-Poisson}.
Moreover, either $u_0$ is periodic (possibly identically equal to zero) or, $u_0>0$ and $u'_0>0$ on $(0,+\infty)$. 	}
 \item [(iv)] if $f \in C^1([0,+\infty))$ with $f(0)<0$ and $f'(t) \geq \alpha >0$ for any $t>0$, then $u$ is one-dimensional and periodic.
\end{enumerate}
}

\textbf{Theorem D} \textit{
Assume $N \geq 2$. The only classical solution to
\begin{equation*}
\left\{
	\begin{array}{ccl}
		-\Delta u= u-1& \text{in}&  \R^{N}_{+},\\
		u\geq 0 & \text{in} & \R^{N}_{+},\\
		u=0 & \text{on} & \partial \R^{N}_{+}, 
	\end{array}
	\right.
\end{equation*}
is the function 
$$u(x) = 1- \cos(x_N)  \quad \text{for any } x \in  \R^{N}_{+}.$$
}

The case $N=2$ in Theorem D is due to \cite{fs1} (just apply item (iv) of Theorem C with $f(u) = u-1$), while the case $N \geq 3$ is proven in \cite{CEGM}.

It is interesting to note that the one-dimensional symmetry results in the four theorems above were obtained without any additional assumption on the non-negative solution $u$. 
It therefore seems natural to seek to broaden the class of functions $f$ guaranteeing that the only classical solutions, bounded or not, to \eqref{Poisson} are necessarily one-dimensional. 
On the other hand, such a one-dimensional symmetry result is not true for every locally Lipschitz-continuous function $f$, as shown by the function $u(x) = x_N e^{-x_1}$, which solves \eqref{Poisson} with $f(u) = - u$ and $ N \geq 2$. Hence, the following two research lines naturally emerge : 
\begin{enumerate}
    \item [(1)] provide general natural conditions on $f$ guaranteeing that the only classical solutions, bounded or not, to \eqref{Poisson} are necessarily one-dimensional;
    \item [(2)] supply sufficiently general assumptions on classical solutions to \eqref{Poisson} ensuring their one-dimensional symmetry  (either independently of the considered function $f \in Lip_{\text{loc}}([0,+\infty))$, or for an appropriate subclass of locally continuous Lipschitz functions).
\end{enumerate}

The new results presented in this article are devoted to both lines of research. 
Before stating our results, let us review those that exist in the literature and that relate to research line (2) above.

With regard to the second line of research, various additional hypotheses on $u$ have been considered in the literature in order to establish their one-dimensional symmetry. They mainly concern the growth conditions of $u$ or its qualitative or spectral properties, such as monotonicity, stability and/or the Morse index. These conditions are natural, as they are generally satisfied by the solutions that appear in applications and by those constructed using variational or topological methods. In this regard, we present below the main contributions found in the literature. We shall describe them in ascending order of generality. Let us begin with the case of bounded solutions.

\textbf{Theorem E} (\cite{BCN}) \textit{
Assume $N=2$ or $3$. Let $u \in C^{2}(\overline{\R^{N}_{+}})$ be a positive and bounded solution to \eqref{Poisson} where 
$f \in Lip_{\text{loc}}([0,+\infty))$. \\
If $N = 2$, $u$ is one-dimensional and strictly increasing. If $N = 3$ the same conclusion holds, if one assumes in addition that $f(0) \geq 0$ and that $f$ is 
$C^1$. 
}

More recently, by combining Theorem C above and a geometric tool introduced in \cite{FarValarma}, the authors of \cite{fs1} improved the result of Theorem E in the two-dimensional case. Indeed, the authors of \cite{fs1} prove the one-dimensional symmetry for any non-negative solution and assuming only that $u$ has a bounded gradient. More precisely,

\textbf{Theorem F} (\cite{fs1}) \textit{
Let $u \in C^{2}(\overline{\R^{2}_{+}})$ be a solution to \eqref{Poisson} where $f \in Lip_{\text{loc}}([0,+\infty))$. \\
If $u$ has bounded gradient, then it is one-dimensional.
}

Recall that, bounded solutions of \eqref{Poisson} have bounded gradient, but the converse is not true, as shown by the (one-dimensional) harmonic function $u(x) = x_N$. Furthermore, Theorem F also covers the case of one-dimensional periodic solutions which, given the homogeneous Dirichlet condition, must necessarily vanish somewhere in $\R^{2}_{+}$.\footnote{This fact is well illustrated by the explicit example discussed in Theorem D.} 

To conclude on the case of bounded solutions of \eqref{Poisson} we mention the two following results. The first one deals with convex functions $f$ with $f(0)=0$, while the second one applies, among other things, to any non-negative function $f$.

\textbf{Theorem G} (\cite{CSlin}) \textit{
Assume $ N \geq 2$ and let $f \in C^1([0,+\infty))\cap C^2((0,+\infty))$ with $f(0) = 0$ and $f''\geq 0$. 
The only bounded classical solution to \eqref{Poisson} is $u \equiv 0$.
}

As an immediate consequence of Theorem G one has 

\textbf{Corollary A} (\cite{CSlin}) \textit{
Assume $N \geq 2$ and let $p>1$. The only bounded classical solution to
\begin{equation*}
\left\{
	\begin{array}{ccl}
		-\Delta u= u^p& \text{in}&  \R^{N}_{+},\\
		u \geq 0 & \text{in} & \R^{N}_{+},\\
		u=0 & \text{on} & \partial \R^{N}_{+},
	\end{array}
	\right.
\end{equation*}	
is $ u \equiv 0$. 
}

Corollary A was first proven for $1 < p < (N+1)/(N-3)^+$ 
in \cite{dan0}, and then for any $1 < p < p_{JL}(N)$ in 
\cite{far}, where $p_{JL}$ is the Joseph-Lundgren stability exponent given by $p_{JL}(N)=\frac{(N-2)^2 -4N + 8 \sqrt{N-1}}{(N-2)(N-10)^+}$.

\smallskip

\textbf{Theorem H} (\cite{far-dup}) \textit{
Let $u\in C^2({\overline {\R^N_+}})$ be a bounded solution of \eqref{Poisson}. Assume $f \in C^1([0,+\infty))$ and 
\begin{enumerate}
\item either $ f(t) \geq 0 $ for $ t \geq 0,$ 
\item or there exists $z>0$ such that $ f(t) \geq 0$ for $ t \in [0,z]$ and $ f(t) \leq 0$ for $t \geq z$. 
\end{enumerate}
If $ 2 \leq N\le 11$, then either $u \equiv 0$, or $u$ is positive, one-dimensional and strictly increasing. 
}

\smallskip

To continue, recall that a solution to \eqref{Poisson} that is either strictly monotone in the direction $x_N$ or a local minimizer is automatically \textit{stable}, which means that it additionally satisfies
\begin{equation}\label{def stable}
    \int_{\R^N_+}|\nabla\varphi|^2 - f'(u) \varphi^2 \geq 0, 
\end{equation}
for all $\varphi\in C^1_c(\R^N_+)$ (here we have supposed $f \in C^1$). \\
Also, we shall say that $u$ is \textit{stable outside a compact} $\mathcal{K} \subset \R^N_+$ if the integral inequality \eqref{def stable} holds for any $\varphi\in C^1_c(\R^N_+ \setminus \mathcal{K})$. Clearly, the stability implies the stability outside a compact set (actually, any finite Morse index solution to \eqref{Poisson} is stable outside a compact set of $\R^N_+$). 

As for these classes of solutions to \eqref{Poisson}, we have

\textbf{Theorem I} (\cite{dup-sir-sou}) \textit{
Assume $N \geq 2$. Let $f \in C^1([0,+\infty))\cap C^2((0,+\infty))$ be non-negative and convex, with $f(0) = 0$ and $f \not  \equiv 0$. \\
Let $u$ be a classical solution of \eqref{Poisson} which is monotone in the $x_N$ direction, i.e. such that $\frac{\partial u}{\partial x_N} \geq 0$ on $\R^{N}_{+}$. Then, $ u \equiv 0$. \\
In particular, $ u \equiv 0$ is the only classical solution to \eqref{Poisson} which is bounded on finite strips.\footnote{ \, When 
$f(0) \geq 0$, any positive solution to \eqref{Poisson}, which is also bounded on finite strips, is strictly monotone in the $x_N-$direction (see Theorem 3.1 of \cite{far1}).}
}

As an immediate consequence of Theorem I one has 

\textbf{Corollary B} (\cite{dup-sir-sou}) \textit{
Assume $N \geq 2$ and let $p>1$. Let $u$ be a classical solution of
\begin{equation}\label{Lane-Emden-eq}
\left\{
	\begin{array}{ccl}
		-\Delta u= u^p& \text{in}&  \R^{N}_{+},\\
		u \geq 0 & \text{in} & \R^{N}_{+},\\
		u=0 & \text{on} & \partial \R^{N}_{+},
	\end{array}
	\right.
\end{equation}	
which is monotone in the $x_N$ direction. Then, $ u \equiv 0$. \\
In particular, $ u \equiv 0$ is the only classical solution to \eqref{Lane-Emden-eq} which is bounded on finite strips.
}

Corollary B recovers and improves upon Corollary A above.

Very recently, Corollary B  has been improved in the work \cite{DFTroy}. Specifically, 

\textbf{Theorem J} (\cite{DFTroy}) \textit{
Assume $N \geq 2$ and let $p>1$. The only classical solution to
\begin{equation*}
\left\{
	\begin{array}{ccl}
		-\Delta u= \vert u \vert^{p-1}u& \text{in}&  \R^{N}_{+},\\
		u=0 & \text{on} & \partial \R^{N}_{+},
	\end{array}
	\right.
\end{equation*}	
which is stable outside a compact set of $\R^{N}_{+}$ is $ u \equiv 0$. 
}

Note that Theorem J holds for sign-changing solutions.  It was first proven in \cite{far} for all $N \leq 10$ and, if $N >10$, for $1<p<p_{JL}(N)=\frac{(N-2)^2 -4N + 8 \sqrt{N-1}}{(N-2)(N-10)}$, the Joseph-Lundgren stability exponent.

We are now ready to present our results. The first theorem provides a contribution to the first line of research above. It recovers and improves upon the result stated in item (iii) of the previous Theorem C. 
It only requires that $f$ behaves at least linearly at infinity, i.e., that $f$ satisfies \eqref{f-superlinear-N=2} below. 

\begin{thm}\label{symmetry-N=2}
		Let $u \in C^{2}(\overline{\R^{2}_{+}})$ be a solution of \eqref{Poisson} where $f \in Lip_{\text{loc}}([0,+\infty))$ satisfies
	\begin{equation}\label{f-superlinear-N=2}
		\liminf\limits_{t \to + \infty} \frac{f(t)}{t}>0.  \\
	\end{equation} 
	Then there exists a bounded function $u_{0}:[0,+\infty) \to [0,+\infty)$ such that
	\begin{equation*}
		u(x)=u_{0}(x_2)  \quad \text{for any } x \in \R^{2}_{+}.
	\end{equation*}
Moreover, either $u_0$ is periodic (possibly identically equal to zero) or, $u_0>0$ and $u'_0>0$ on $(0,+\infty)$. 	
\end{thm}

Next we state our results concerning the second line of research above. Let us begin with the case of solutions monotone 
in the $x_N-$direction. 

\begin{thm}\label{symmetry-half-space-u-increase}
	Assume $2 \leq N \leq 9$ and let $u \in C^{2}(\overline{\R^{N}_{+}})$ be a solution of \eqref{Poisson} such that $\frac{\partial u}{\partial x_{N}} \geq 0$ on $\R^{N}_{+}$. 
	Assume that $f \in C^{1}([0,+\infty))$ is a non-negative function satisfying
	\begin{equation}\label{superlinear-u-monotone}
		\liminf\limits_{t \to + \infty} \frac{f(t)}{t}>0.  
	\end{equation} \\
Then, either $u \equiv 0$, or there exists a bounded function $u_{0}:[0,+\infty) \to [0,+\infty)$ such that
	\begin{equation*}
		u(x)=u_{0}(x_{N})  \qquad \text{for any } x \in \R^{N}_{+},
	\end{equation*}
and $u_0>0$, $u'_0>0$ on $(0,+\infty)$. 

Moreover, if $f(t)>0$ for any $t>0$, then necessarily $u\equiv 0$ in $\R^{N}_{+}$ and $ f(0)=0$.
\end{thm}

Although subject to a dimensional restriction, the previous result extends Theorem I to all non-negative functions $f$ under the sole assumption that $f$ behave at least linearly at infinity.\footnote{\, Note that, the convex function $f$ in Theorem I necessarily satisfies the assumption \eqref{superlinear-u-monotone}.}  It would be interesting to know whether this result still holds for 
$N>9$.

If we further assume that $f$ is convex or that $f$ has at most polynomial growth, the conclusions of Theorem \ref{symmetry-half-space-u-increase} remain valid up to dimension $N=10$. Specifically,

\begin{thm}\label{symmetry-half-space-u-increase2}
	Assume $2 \leq N \leq 10$ and let $u \in C^{2}(\overline{\R^{N}_{+}})$ be a solution of \eqref{Poisson} such that $\frac{\partial u}{\partial x_{N}} \geq 0$ on $\R^{N}_{+}$. 
	Assume that $f \in C^{1}([0,+\infty))$ is a non-negative, convex function satisfying
	\begin{equation}\label{superlinear-u-monotone2}
		\liminf\limits_{t \to + \infty} \frac{f(t)}{t}>0.  
	\end{equation} \\	
Then, either $u \equiv 0$, or there exists a bounded function $u_{0}:[0,+\infty) \to [0,+\infty)$ such that
	\begin{equation*}
		u(x)=u_{0}(x_{N})  \qquad \text{for any } x \in \R^{N}_{+},
	\end{equation*}
and $u_0>0$, $u'_0>0$ on $(0,+\infty)$. 

Moreover, if $f(0)=0$, then necessarily $u\equiv 0$ in $\R^{N}_{+}$.
\end{thm}

Unlike Theorem I, the function $f$ in Theorem \ref{symmetry-half-space-u-increase2} is not necessarily non-decreasing. For instance, Theorem \ref{symmetry-half-space-u-increase2} applies to $f(t) = \vert t-1\vert^p$, $p>1$.

We also note that, if we drop the assumption \eqref{superlinear-u-monotone2}, then the full conclusion of Theorem \ref{symmetry-half-space-u-increase2} would no longer be valid. Indeed, when $f(t) = e^{-2t}$, which is positive and convex, the \textit{unbounded} and monotone function $u(x)=\ln{(1 + x_N)}$ solves  \eqref{Poisson} for any $ N\geq2.$ \\
This observation also applies to Theorem \ref{symmetry-half-space-u-increase} above and to the following theorem. 

\medskip

\begin{thm}\label{symmetry-half-space-u-increase3}
	Assume $2 \leq N \leq 10$ and let 
	$u \in C^{2}(\overline{\R^{N}_{+}})$ be a solution of \eqref{Poisson} such that $\frac{\partial u}{\partial x_{N}} \geq 0$ on $\R^{N}_{+}$. 
	Assume that $f \in C^{1}([0,+\infty))$ is a non-negative function satisfying 
\begin{equation}\label{superaff-u-croiss-contr2}
    \liminf\limits_{t \to + \infty} \frac{f(t)}{t} >0,
    \qquad \limsup\limits_{t \to + \infty} \frac{f(t)}{t^\theta}< +\infty,
	\end{equation}
	for some $\theta>1$. \\	
Then, either $u \equiv 0$, or there exists a bounded function $u_{0}:[0,+\infty) \to [0,+\infty)$ such that
	\begin{equation*}
		u(x)=u_{0}(x_{N})  \qquad \text{for any } x \in \R^{N}_{+},
	\end{equation*}
and $u_0>0$, $u'_0>0$ on $(0,+\infty)$. 

Moreover, if $f(t)>0$ for any $t>0$, then necessarily $u\equiv 0$ in $\R^{N}_{+}$ and $ f(0)=0$.
\end{thm}

As an immediate consequence of previous three theorems we can treat the case of solutions that are bounded on finite strips. 

\begin{cor}\label{symmetry-half-space-u-bound-strip}
Let $N$, $f$ be as in the statement of Theorem \ref{symmetry-half-space-u-increase} (resp. Theorem \ref{symmetry-half-space-u-increase2} or Theorem \ref{symmetry-half-space-u-increase3}) and let $u$ be a classical solution of \eqref{Poisson} which is bounded on finite strips. \\	Then, the conclusions of Theorem \ref{symmetry-half-space-u-increase} (resp. Theorem \ref{symmetry-half-space-u-increase2} or Theorem \ref{symmetry-half-space-u-increase3}) still hold true.
\end{cor}

The viewpoint we adopted to prove the previous theorems can also be applied to both stable solutions and stable solutions outside a compact set. Below we present the Liouville-type theorems corresponding to these two classes of solutions.

\begin{thm}\label{symmetry-half-space-superlinear}
	Assume $2 \leq N \leq 4$ and let $u \in C^{2}(\overline{\R^{N}_{+}})$ be a stable solution of \eqref{Poisson}
	where $f \in C^{1}([0,+\infty))$ is a non-negative and non-decreasing function satisfying
	\begin{equation}\label{superlinear1}
	\lim\limits_{t \to + \infty} \frac{f(t)}{t}=+\infty.
	\end{equation}
Then, $u\equiv 0$ in $\R^{N}_{+}$ and $ f(0) = f'(0) =0$.
\end{thm}

\medskip

In the next result we replace the superlinearity condition \eqref{superlinear1} on $f$ with a growth assumption on the solution $u$. 

\medskip

\begin{thm}\label{symmetry-half-space-N=3,4-dirichlet}
	Assume $2 \leq N \leq 4$ and let 
	$u \in C^{2}(\overline{\R^{N}_{+}})$ be a stable solution of \eqref{Poisson}
	where $f \in C^{1}([0,+\infty))$ is a  non-negative, non-decreasing function with $ f \not \equiv 0$. If 
	\begin{equation}\label{grand-O-half-space}
		u(x)=O(|x|^{\frac{4-N}{2}} \ln^{1/2} |x|) \qquad \text{as} \quad |x| \to + \infty,
	\end{equation}
	then, $u\equiv 0$ in $\R^{N}_{+}$ and $ f(0) = f'(0) =0$.
\end{thm}

\medskip

The following two results deal with solutions which are stable outside a compact set of $\R^{N}_{+}$. 

\medskip

\begin{thm}\label{symmetry-half-space-f-convex-increasing}
	Assume $2 \leq N \leq 9$ and let $u \in C^{2}(\overline{\R^{N}_{+}})$ be a solution of \eqref{Poisson} which is stable outside a compact set of $\R^{N}_{+}$. Let $f \in C^{1}([0,+\infty))$ be a non-negative, convex function with $f(0)=0$ and 
	$f \not \equiv 0$. \\
Then, $u\equiv 0$ in $\R^{N}_{+}$ and $f'(0) =0$.
\end{thm}

\medskip

\begin{thm}\label{symmetry-half-space-f-convex-increasing-croiss-contr}
	Assume $2 \leq N \leq 10$ and let 
	$u \in C^{2}(\overline{\R^{N}_{+}})$ be a solution of \eqref{Poisson} which is stable outside a compact set of $\R^{N}_{+},$ where $f \in C^{1}([0,+\infty))$ is a non-negative and non-decreasing function satisfying 
\begin{equation}\label{superaff-u-croiss-contr}
    \lim\limits_{t \to + \infty} \frac{f(t)}{t}=+\infty, 
    \qquad \limsup\limits_{t \to + \infty} \frac{f(t)}{t^\theta}< +\infty,
	\end{equation}
	for some $\theta>1$. \\	
Then, $u\equiv 0$ in $\R^{N}_{+}$ and $f(0) = f'(0) =0$.
\end{thm}

In the case of monotone solutions, the flexibility of our approach allows us to extend some of the previous results to the case of differential inequalities, even without imposing any boundary conditions. More precisely, we have 

\medskip

\begin{thm}\label{sym-diff-ineq-1}
	Assume $N = 2,3 $ and let $u \in C^2(\R^{N}_{+})$ be a solution of 
\begin{equation}\label{ineq-Poisson-gen-1}	
	\left\{
	\begin{array}{ccl}
		-\Delta u \geq f(u) & \text{in}&  \R^{N}_{+},\\
		u \geq 0 & \text{in} & \R^{N}_{+},\\
		\frac{\partial u}{\partial x_{N}} \geq 0 & \text{in} & \R^{N}_{+},
	\end{array}
	\right.
\end{equation} 		
where $f \in C^0([0,+\infty))$ satisfies $f(t)>0$ for any $t>0$ and 
	\begin{equation*}
		\liminf\limits_{t \to + \infty} \frac{f(t)}{t}>0.  \\
	\end{equation*} 
	Then, $u\equiv 0$ in $\R^{N}_{+}$ and $f(0)=0$.
\end{thm}

In particular, for $N=2,3$ and for any $p \geq 1$, the only non-negative and monotone solution to $-\Delta u \geq u^p $ in $\R^{N}_{+}$ is the function identically equal to zero. 

We also note that the result above is sharp. Indeed, the function 
$u(x',x_N)= c(1 + \vert x' \vert^2)^{-\frac{1}{p-1}}$ satisfies $-\Delta u \geq u^p $ and $\frac{\partial u}{\partial x_{N}} \equiv 0$ in $\R^{N}_{+}$, when $N \geq 4$, $ p>\frac{N-1}{N-3}$ and $c>0$ is small enough. Nevertheless, we have the following sharp Liouville-type result. 

\begin{thm}\label{sym-diff-ineq-2}
Assume $N \geq 2$, let $p\geq 1$ be a real number and let $u \in C^2(\R^{N}_{+})$ be a solution of
\begin{equation}\label{ineq-Lane-Emden-gen-1}	
	\left\{
	\begin{array}{ccl}
		-\Delta u \geq u^p & \text{in}&  \R^{N}_{+},\\
		u \geq 0 & \text{in} & \R^{N}_{+},\\
		\frac{\partial u}{\partial x_{N}} \geq 0 & \text{in} & \R^{N}_{+}.
	\end{array}
	\right.
\end{equation} 	
If $1 \leq p \leq \bar p(N)$, where 
\begin{equation}\label{p_c(N)}	
\bar p(N):= \left\{
			\begin{array}{crl}
      		+\infty & \text{if}& N=2,3,\\
        \frac{N-1}{N-3} & \text{if} & N \geq 4,
			\end{array}
				\right.
\end{equation}
then, $u\equiv 0$ in $\R^{N}_{+}$.

\end{thm}

If we restrict ourselves to solutions of the semilinear Poisson equation on $\R^{N}_{+}$, but still without imposing boundary conditions, the classification results above can be extended to dimension $ N \geq 4$. This is the content of the following two theorems.   

\begin{thm}\label{sym-diff-ineq-3}
	Assume $2 \leq N \leq 9$ and let $u \in C^2(\R^{N}_{+})$ be a solution of 
\begin{equation}\label{ineq-Poisson-gen-2}	
	\left\{
	\begin{array}{ccl}
		-\Delta u = f(u) & \text{in}&  \R^{N}_{+},\\
		u \geq 0 & \text{in} & \R^{N}_{+},\\
		\frac{\partial u}{\partial x_{N}} \geq 0 & \text{in} & \R^{N}_{+},
	\end{array}
	\right.
\end{equation} 		
where $f \in C^1([0,+\infty))$ satisfies $f(t)>0$ for any $t>0$ and 
	\begin{equation}\label{f-superlinear-N=2-ineq}
		\liminf\limits_{t \to + \infty} \frac{f(t)}{t}>0.  \\
	\end{equation} 
Then, either $u\equiv 0$ in $\R^{N}_{+}$, or there exists a function $u_{0}:\R^{N-1} \to (0,+\infty)$ such that 
\begin{equation*}
		u(x)=u_{0}(x_1,..., x_{N-1}) \qquad \text{for any } x \in \R^{N}_{+}. 
\end{equation*}
If we further assume \eqref{superaff-u-croiss-contr2}, the above conclusion holds true even in dimension $N=10$.
\end{thm}

In the case of the Lane-Emden equation $-\Delta u = u^p$, the above result can be improved as follows 

\begin{thm}\label{sym-diff-ineq-4} Assume $N \geq 2$ and let $u \in C^2(\R^{N}_{+})$ be a solution of
\begin{equation}\label{ineq-Lane-Emden-gen-2}	
	\left\{
	\begin{array}{ccl}
		-\Delta u = u^p & \text{in}&  \R^{N}_{+},\\
		u \geq 0 & \text{in} & \R^{N}_{+},\\
		\frac{\partial u}{\partial x_{N}} \geq 0 & \text{in} & \R^{N}_{+}. 
	\end{array}
	\right.
\end{equation} 	
If $N=2,3$, then $u\equiv 0$ in $\R^{N}_{+}$; while for $N \geq 4$ we have 

\begin{enumerate}
    \item [(i)] if $1 \leq p < \frac{N+1}{N-3}$, then $u\equiv 0$ in $\R^{N}_{+}$;

\bigskip

  \item [(ii)] if $p= \frac{N+1}{N-3}$, then either $u\equiv 0$ in $\R^{N}_{+}$ or \\
\begin{equation}\label{Talenti-Aubin-funct}
		u(x)= (a + b \vert x'-x'_0 \vert^2)^{-\frac{N-3}{2}} \qquad \text{for any } x = (x',x_N)\in \R^{N}_{+}, 
\end{equation}
for $a,b>0$ satisfying $1=ab(N-1)(N-3)$ and any $x'_0 \in \R^{N-1}$; 

\bigskip

\item [(iii)] if $\frac{N+1}{N-3} < p < p_{JL}(N-1)$, where $p_{JL}(N)=\frac{(N-2)^2 -4N + 8 \sqrt{N-1}}{(N-2)(N-10)^+}$ is the Joseph-Lundgren stability exponent. \\Then 
either $u\equiv 0$ in $\R^{N}_{+}$, or there exists a function $u_{0}:\R^{N-1} \to (0,+\infty)$ such that 
\begin{equation*}
		u(x)=u_{0}(x') \qquad \text{for any } x = (x',x_N)\in \R^{N}_{+},
\end{equation*}
\end{enumerate}
\end{thm}

When $N\geq 4$, 
for any $\frac{N+1}{N-3} < p < p_{JL}(N-1)$, there exists a positive solution to \eqref{ineq-Lane-Emden-gen-2} depending only on the first $N-1$ variables.\footnote{\, 
For any $N \geq 3$  and $p > \frac{N+2}{N-2}$ positive radial solutions to $-\Delta u =u^p$ in $\R^N$ always exist (see for instance  
\cite{GS2},\cite{JL}). If $u_0$ is such a function in dimension $N-1$, then the function $u(x)  = u_0(x')$, $x = (x',x_N) \in \R^{N}_{+}$, provides the desired example.} 
Also note that, the positive solutions in items $(ii)$ and $(iii)$ are unstable (see Theorem 5 and Remark 4 of \cite{far}). Finally observe that, for $ N \leq 11$, the result above provides the complete classification of the solutions to \eqref{ineq-Lane-Emden-gen-2} for any $p \geq 1$.

\section{Proofs of results in the case of the semilinear Poisson equation}

\begin{proof}[Proof of Theorem \ref{symmetry-N=2}]
By the results in \cite{fs1,fs2} (see Theorem C in the Introduction), either $u$ is one-dimensional and periodic (possibly identically equal to zero) and we are done, 
or $ u > 0 $ and $\frac{\partial u}{\partial x_2} > 0$ in $\R^2_{+}$. 
In the latter case, by Lemma \ref{lemma-profile-infty} the profile of $u$ at infinity defined by 
$$
\overline{u}(x_1) := \lim_{x_2 \to + \infty} u(x_1,x_2), \qquad \text{for any } x_1 \in \R,
$$
is a non-negative member of $L^1_{loc}(\R)$ such that $f(\overline{u}) \in L^1_{loc}(\R)$ and 
\begin{equation}\label{1D-pbl}
- \overline{u}{''} \geq  f(\overline{u}) \qquad \text{in} \quad  \mathcal{D}^{'}(\R).\\
\end{equation}
Furthermore, the sequence of functions 
$(u_n)_{n \geq 1}$ defined by 
$$
u_n(x) := \tilde u(x_1,x_2 + n), \qquad \text{for any } x \in \R^{2},
$$
converges to $\overline{u}$ in $L^1_{loc}(\R^{2})$. 

From \eqref{1D-pbl} and \eqref{f-superlinear-N=2} we see that 
\begin{equation}\label{1D-ineg}
- \overline{u}{''} \geq  A \overline{u} -B  \qquad \text{in} \quad  \mathcal{D}^{'}(\R), \\
\end{equation}
for some constants $A>0$ and $B \geq 0$. Then the function $v(x_1) = - \overline{u}(x_1) + \frac{B}{2} x_1^2$ belongs to $L^1_{loc}(\R)$ and solves 
\begin{equation}\label{1D-ineg-bis}
- v{''} \leq  0 \qquad \text{in} \quad  \mathcal{D}^{'}(\R), \\
\end{equation}
hence, $v$ is subharmonic in the sense of the Classical Potential Theory (see for instance Theorem 3.2.1 in \cite{DiP-Val} or Theorem 3.2.11 in \cite{Horm}). Then $v$ is convex on $\R$ (see for instance Lemma 3.2.9 in  \cite{DiP-Val}) thus continuous on $\R$ and so is $\overline{u}$. 
We observe that, for any compact set $K \subset \R^{2}$ we have $ 0 \leq u_n \leq \overline{u} \leq \|\overline{u}\|_{L^{\infty}(K)} \textbf{1}_K$, since $\frac{\partial u}{\partial x_2} > 0 \in \R^2_{+}$, so $\vert f(\overline{u}) -f(u_n) \vert \leq L \vert \overline{u}-u_n \vert$ on $K$, where $L$ is the Lipschitz constant of $f$ in the compact interval $[0,\|\overline{u}\|_{L^{\infty}(K)}]$. Then, $f(u_n) \longrightarrow f(\overline{u})$ in $L^1_{loc}(\R^{2})$ so that we can pass to the limit in the equation satisfied by $u_n$ to get that
\begin{equation}\label{1D-eg-tris}
- \overline{u}{''} = f(\overline{u})  \qquad \text{in} \quad  \mathcal{D}^{'}(\R). \\
\end{equation}
Since $f(\overline{u})$ is continuous, from \eqref{1D-eg-tris} we infer that $\overline{u} \in C^2(\R)$ is, indeed, a classical solution. From the latter property and \eqref{1D-ineg} we see that $\overline{u}$ is bounded\footnote{\, Indeed, $\overline{u}$ satisfies $\frac{(\overline{u}'(t))^2}{2} + F(\overline{u}(t)) = \frac{(\overline{u}'(0))^2}{2} +  F(\overline{u}(0))$, for any $t \in \R$, where $F(u) := \int_0^u f(s) ds. $ Hence, $ A \frac{(\overline{u}(t))^2}{2} - B \overline{u}(t) \leq \frac{(\overline{u}'(0))^2}{2} +  F(\overline{u}(0))$. The boundedness of $\overline{u}$ then follows, since $A>0$. } on $\R$ which, in turns, implies that also $u$ is bounded on $\R^{2}_+.$ Then, from Theorem E (or Theorem F) in the Introduction, it follows that $u$ is one-dimensional. \end{proof}

\medskip

\begin{proof}[Proof of Theorem \ref{symmetry-half-space-u-increase}]
By the strong maximum principle, either $ u \equiv 0 $ (and we are done) or $u > 0 $ and $\frac{\partial u}{\partial x_{N}} > 0$ in $\R^{N}_{+}$. The latter implies that $u$ is stable on $\R^{N}_{+}$ and so Theorem 1.2 in \cite{crfs} tell us that
\begin{equation}\label{estimate-alpha-holder-1}
    \|u\|_{L^{\infty}({B(x,\frac{R_0}{2}}))} \leq C \|u\|_{L^{1}(B(x,R_0))} \qquad \text{for any} \quad B(x,R_0) \subset \subset \R^{N}_{+},
\end{equation}
where  $C=C(N,R_0)>0$ is constant.\\
The latter and the assumption \eqref{superlinear-u-monotone} enable us to apply item (i) of Lemma \ref{lemma-bound-L-1} to deduce that $u$ is bounded on the half-space $\{x \in \R^{N}_{+} : x_N > 2R_0 \}$. The monotonicity of $u$ then implies that $u$ is bounded on $\R^{N}_{+}$. The desired conclusion the follows by invoking Theorem 4 in \cite{far-dup} (see Theorem H in the Introduction). \\
Let us now move on to the last statement of the theorem and show that $u \equiv 0$ if $f(t)>0$ for any $t>0$. If not, then $u=u_0(x_N)$, where $u_{0}:[0,+\infty) \to [0, ,+\infty)$ is a bounded, positive and strictly increasing solution to the ODE $-u'' =f(u)$ on $(0,+\infty)$. By the first integral of the ode we have that also $u'$ is bounded and, from this information, we immediately infer that $f(\sup u_0)=0$. A contradiction.
\end{proof}

\medskip

\begin{proof}[Proof of Theorem \ref{symmetry-half-space-u-increase2}]
As in the proof of Theorem \ref{symmetry-half-space-u-increase}, either $ u \equiv 0 $ (and we are done) or $u > 0,$ 
$\frac{\partial u}{\partial x_{N}} > 0 $ in $\R^{N}_{+}$ and $u$ is stable on $\R^{N}_{+}$. Also, from \eqref{superlinear-u-monotone2}, we have
\begin{equation}\label{f-superlinear-2}
	f(t)\geq At-B \qquad \text{ for any } t>0,
\end{equation}
for some constants $A>0$ and $B \geq 0$. \\
By Lemma \ref{lemma-profile-infty} the profile of $u$ at infinity defined by 
$$
\overline{u}(x') := \lim_{x_N \to + \infty} u(x',x_N), \qquad \text{for any } x' \in \R^{N-1},
$$
is a non-negative member of $L^1_{loc}(\R^{N-1})$ and  $f(\overline{u}) \in L^1_{loc}(\R^{N-1})$. Furthermore, the sequence of non-negative functions $(u_n)_{n \geq 1}$ defined by 
$$
u_n(x) := \tilde u(x',x_N + n), \qquad \text{for any } x \in \R^N,
$$
converges to $\overline{u}$ in $L^1_{loc}(\R^{N})$. Also observe that $u_n \in C^0(\R^N) \cap H^1_{loc}(\R^{N})$ for any $n \geq 1.$

Let us prove that $\overline{u}\in H^1_{loc}(\R^{N-1})$ and that $\overline{u}$ is a stable solution to $ -\Delta \overline{u} = f(\overline{u})$ in $\mathcal{D}^{'}(\R^{N-1}).$ 

Under the assumptions on $f$, there exists $t_0 >0$ such that $f'(t)>0$ on $(t_0, \infty)$. Then
$$
\vert f(u_n) \vert  = f(u_n) = f(u_n) \textbf{1}_{\{u_n \leq t_0\}} +  f(u_n) \textbf{1}_{\{u_n >t_0\}}  \leq 
\Vert f \Vert_{L^{\infty}([0,t_0])} + f(\overline{u})\in L^1_{loc}(\R^{N}),
$$
and so 
\begin{equation}\label{cv-f(u-n)}
f(u_n) \longrightarrow f(\overline{u}) \qquad \text{in} \quad L^1_{loc}(\R^{N}),
\end{equation}
by the dominated convergence theorem. 

Since $ -\Delta u_n = f(u_n)$ in $\R^{N}_+$ and $u_n \longrightarrow \overline{u} $ in $L^1_{loc}(\R^{N})$, $f(u_n) \longrightarrow f(\overline{u})$ in $L^1_{loc}(\R^{N})$, we get that $\overline{u}$ solves $ -\Delta \overline{u} = f(\overline{u})$ in $\mathcal{D}^{'}(\R^{N})$. The latter implies that $ -\Delta \overline{u} = f(\overline{u})$ in $\mathcal{D}^{'}(\R^{N-1}),$ since $\overline{u}$ depends only on the first $N-1$ variables. 

To prove that $\overline{u} \in H^1_{loc}(\R^{N-1})$ we recall that $u$ is stable on $\R^{N}_+.$ 
Therefore, Proposition 2.5 in \cite{crfs} tell us that 
    \begin{equation}\label{estimate-Sobolev-1}
		\int_{B(z,R)}|\nabla u |^{2} \leq 
		C' R^{N-2}\Big( \frac{1}{R^{N}} \int_{B(z,2R)} u \Big)^{2} \qquad \text{for any } B(z,2R) \subset \subset \R^{N}_{+},
	\end{equation}
where $C'>0$ is a dimensional constant.

Set $e_N =(0,...,0,1) \in \R^N$, pick any $ r \geq R_0$, where $R_0$ is given by item (i) of Lemma \ref{lemma-bound-L-1} and consider any ball $B(z,r) \subset \R^{N}$.  There exists an integer $\bar n$, depending only on $z$ and $r$, such that 
$ B(z+n e_N, 4r) \subset \subset \R^{N}_{+}$ for any $ n \geq \bar n$.  Hence, by \eqref{estimate-Sobolev-1} (with $R=r$) we deduce that 
$$
\int_{B(z, r)} |\nabla u_n |^{2} =\int_{B(z+n e_N, r)} |\nabla u |^{2} 
\leq C' r^{N-2}\Big( \frac{1}{r^{N}} \int_{B(z+n e_N,2r)} u \Big)^{2} \qquad \text{for any } n \geq \bar n.
$$
Now, for any $n \geq \bar n$, an application of Lemma \ref{lemma-bound-L-1} (with $R=2r$) yields
\begin{equation}\label{estimate-Sobolev-2-2}
\int_{B(z, r)} |\nabla u_n |^{2} \leq C' r^{N-2}\Big( \frac{1}{r^{N}} \int_{B(z+n e_N,2r)} u \Big)^{2}  \leq C' r^{N-2} 
\Big( \frac{C_0 (2r)^N}{r^{N}} \Big)^{2} = C_1(f,N) r^{N-2},
	\end{equation}
where $C_1=C_1(f,N)>0$ is a constant (independent of $n$). 


Moreover, if we denote by $[u_{z,n,r}]$ the average value of $u$ over the open ball $B(z+n e_N,r)$ then, for any $n \geq \bar n$, the Poincaré inequality gives  
    \begin{equation}\label{estimate-Poincaré}
    	\begin{split}
         \int_{B(z, r)} u_n^2  & = \int_{B(z+n e_N, r)} u^2 \leq 2 \int_{B(z+n e_N, r)} \vert u - [u_{z,n,r}] \vert^2 + 2 \int_{B(z+n e_N, r)}  [u_{z,n,r}]^2 \\ 
         & \leq C(N)r^2 \int_{{B(z+n e_N, r)} }|\nabla u |^{2} + C(N)r^N \Big( \frac{1}{r^{N}} \int_{{B(z+n e_N, r)}} u \Big)^{2} \\ 
         & \leq C(N)r^2 \int_{{B(z, r)} }|\nabla u_n |^{2} + C(N)r^N \Big( \frac{1}{r^{N}} \int_{{B(z+n e_N, r)}} u \Big)^{2} \\  
         & \leq C(N)r^2 C_1 r^{N-2} +  C(N)C_0(f,N)^2 r^N = C_2(f,N) r^N,
         \end{split}
	\end{equation}
where in the latter we have used \eqref{estimate-Sobolev-2-2}, as well as Lemma \ref{lemma-bound-L-1} with $R=r$, since we also have that  $B(z+ n e_N, 2r) \subset \subset \R^{N}_+$.
	
By summarizing, we have proved that , 
\begin{equation}\label{estimate-Sobolev-2}
		\int_{B(z,r)}|\nabla u_n |^{2} + u_n^{2} \leq C_3(f,N,r) \qquad \text{for any } n \geq \bar n, 
	\end{equation}
where $C_3>0$ is independent of $n$. Hence, the sequence $(u_n)$ is bounded in $H^1_{loc} (\R^{N})$ and so, we can find a subsequence, still denoted by $(u_n)$, and $v \in H^1_{loc}(\R^{N})$ such that
\begin{equation}\label{cv-faible-u-n}
u_n \longrightarrow v \qquad \text{weakly in} \quad H^1_{loc} (\R^{N}).
\end{equation}
The latter and the fact that $u_n \longrightarrow \overline{u}$ in $L^1_{loc}(\R^{N})$ immediately imply that $ \overline{u} = v \in H^1_{loc}(\R^{N})$. Hence, $\overline{u} \in H^1_{loc}(\R^{N-1})$, since $\overline{u}$ is independent of the variable $x_N$. 

Next we prove that  $\overline{u}$ is a stable solution to $ -\Delta \overline{u} = f(\overline{u})$ in $\mathcal{D}^{'}(\R^{N-1}).$ 
To this purpose, we observe that $u_n$ is strictly increasing in $\R^{N}_{+}$, hence stable, therefore  
\begin{equation}\label{stab-u_n}
\int f'(u_n) \varphi^2 \leq \int |\nabla \varphi|^2 \qquad \text{for any } \varphi \in C^{\infty}_c(\R^{N}_{+}),
\end{equation}
and so
\begin{equation}
\int \vert f'(u_n) \vert \varphi^2 = \int f'(u_n) \varphi^2 + 2 \int [f'(u_n)]^- \varphi^2
\leq \int |\nabla \varphi|^2 + \int 2 [f'(0)]^- \varphi^2
\end{equation}
by the convexity of $f$. Hence, by Fatou's Lemma, we have  
\begin{equation}\label{stab-u_n-bis}
\int \vert f'(\overline{u}) \vert \varphi^2 
\leq \int |\nabla \varphi|^2 + \int 2 [f'(0)]^- \varphi^2 \qquad \text{for any } \varphi \in C^{\infty}_c(\R^{N}_{+}).
\end{equation}
The latter immediately imply that $f'(\overline{u}) \in L^1_{loc}(\R^{N}_+)$, and so $f'(\overline{u})$ belongs to $L^1_{loc}(\R^{N-1})$ as well as to $L^1_{loc}(\R^{N})$ (since $f'(\overline{u})$ only depends on the first $N-1$ variables).  Also
\begin{equation}
\begin{split}
\vert f'(u_n) \vert   &= \vert f'(u_n) \vert \textbf{1}_{\{u_n \leq t_0\}} + \vert f'(u_n) \vert \textbf{1}_{\{u_n >t_0\}} \\
& \leq 
\Vert f' \Vert_{L^{\infty}([0,t_0])} + f'(u_n) \textbf{1}_{\{u_n >t_0\}} \leq 
\Vert f' \Vert_{L^{\infty}([0,t_0])} +  \vert f'(\overline{u}) \vert  \in L^1_{loc}(\R^{N})
\end{split}
\end{equation}
and so, $f'(u_n) \longrightarrow f'(\overline{u})$ in $L^1_{loc}(\R^{N})$ by the dominated convergence theorem. 
Then, passing to the limit in $\eqref{stab-u_n}$ we get 
\begin{equation}\label{quasi-stab-u-lim}
\int f'(\overline{u}) \varphi^2 \leq \int |\nabla \varphi|^2 \qquad \text{for any } \varphi \in C^{\infty}_c(\R^{N}_{+}),
\end{equation}

Now we consider a function $\eta \in C^{\infty}_c(\R_{+})$ such that $supp \, \eta \subset (-1,1)$  and 
$\int_{\R} \eta^2 =1$, and we insert in the stability inequality \eqref{quasi-stab-u-lim} the test function $\varphi_n (x',x_N) := 
\phi (x') \eta_n(x_N)$, where $\phi \in C^{\infty}_c(\R^{N-1})$ and $\eta_n = \frac{1}{n} \eta(\frac{x_N - 4n}{n})$. 
Note that this choice is possible since $\varphi_n$ is supported in $\R^{N}_{+}$ by construction. Hence
\begin{equation}\label{stab-limite}
\begin{split}
&\int_{\R^{N}_{+}} f'(\overline{u}) \phi^2 \eta_n^2 \leq \int_{\R^{N}_{+}} |\nabla \phi|^2 \eta_n^2 + \int_{\R^{N}_{+}} \phi^2 (\eta^{'}_n)^2 \\
& (\int_{\R^{N-1}} f'(\overline{u}) \phi^2) (\int_{\R_+} \eta_n^2) \leq (\int_{\R^{N-1}} |\nabla \phi|^2) (\int_{\R_+}\eta_n^2) + 
(\int_{\R^{N-1}} \phi^2) (\int_{\R_+} (\eta^{'}_n)^2 ) \\
& \int_{\R^{N-1}} f'(\overline{u}) \phi^2 \leq \int_{\R^{N-1}} |\nabla \phi|^2 + 
(\int_{\R^{N-1}} \phi^2) (\int_{\R_+} (\eta^{'}_n)^2 ) \\
\end{split}
\end{equation}
and the last term in the latter tends to zero, when $ n \to \infty$, since $\int_{\R_+} (\eta^{'}_n)^2 = \frac{1}{n^2} 
\int_{\R} (\eta^{'})^2 \leq \frac{2\Vert \eta^{'}\Vert^{2}_{L^{\infty}(\R_{+})}}{n}.$ Passing to the limit into \eqref{stab-limite} we immediately get the stability of  $\overline{u}$. 

By summarizing, we have proved that $\overline{u}\in H^1_{loc}(\R^{N-1})$ is a stable solution to $ -\Delta \overline{u} = f(\overline{u})$ in $\mathcal{D}^{'}(\R^{N-1})$. Since $ N-1 \leq 9$, we can apply Theorem 5 in \cite{far-dup2} to conclude that $\overline{u} \in C^2(\R^{N-1})$, i.e., $\overline{u}$ is a classical stable solution to $ -\Delta \overline{u} = f(\overline{u})$ on $\R^{N-1}$. We can therefore invoke Theorem 1 (or Theorem 2) in \cite{far-dup} to prove that $\overline{u}$ is constant. This result and the monotonicity of $u$ imply that $u$ is bounded on $\R^{N}_{+}$. The desired conclusion then follows by applying Theorem 4 in \cite{far-dup} (see Theorem H in the Introduction).

Let us now move on to the last statement of the theorem and show that $u \equiv 0$ if $f(0)=0$. If not, then $u=u_{0}(x_N)$, where $u_0 : [0,+\infty) \to (0,+\infty)$ is a bounded, positive and strictly increasing function satisfying the ODE $ -u_0^{''} =f(u_0)$ on $[0,+\infty)$. Let us prove that this is impossible. Indeed, by the first integral of the ode we get that also $u'$ is bounded, and so $f(\sup u_{0})=0$. Since $f$ is non-negative and convex with $f(0)=0$, $f$ must be non-decreasing too, so $ 0 \leq f(u_0(t)) \leq f(\sup u_{0})=0 $ for any $t \geq 0$.  Therefore, $ -u_0^{''} = 0$, $u_0>0$ on $(0,+\infty)$ and $u_0(0)=0$, so that $u_{0}(t)=\alpha t $ with $\alpha >0$, contradicting the boundedness of $u_0$.

\end{proof}

\medskip

\begin{proof}[Proof of Theorem \ref{symmetry-half-space-u-increase3}]
Unless we add dummy variables, we may and do assume that $N=10$. 
By the strong maximum principle, either $ u \equiv 0 $ (and we are done) or $u > 0 $ and $\frac{\partial u}{\partial x_{N}} > 0$ in $\R^{N}_{+}$. Thus $u$ is stable on $\R^{N}_{+}$ and so, Theorem 7.1 in \cite{crfs} tell us that
\begin{equation}\label{estimate-any-p}
    \forall \, p>1 \quad \exists C_1=C_1(N,p,R)>0 \quad \text{:} \quad \|u\|_{L^p({B(x,\frac{R}{4}}))} \leq C_1 
    \|u\|_{L^{1}(B(x,R))} \quad 
    \forall \,\, B(x,R) \subset \subset \R^{N}_{+}.
\end{equation}

Then, using the second condition in \eqref{superaff-u-croiss-contr2} we get that 
\begin{equation}\label{estimate-f(u)-any-p>1}
\begin{split}
   &\forall \, p>1 \quad \exists C_2=C_2(N,p,R,f,\theta)>0 \quad \text{such that} \\ &\|f(u)\|_{L^p({B(x,\frac{R}{4}}))} \leq C_2 
   \Big [1 + \|u\|^{\theta}_{L^{p \theta}(B(x,\frac{R}{4}))} \Big]
   \qquad \forall \,\, B(x,R) \subset \subset \R^{N}_{+},
\end{split}   
\end{equation}
where $C_2=C_2(N,p,R_0, f, \theta)>0$ is a constant.\\
Now, the first condition in \eqref{superaff-u-croiss-contr2} and item (i) of Lemma \ref{lemma-bound-L-1} yield
\begin{equation}\label{bound-L-1-stab-out-10}
		\int_{B(x,R_0)}u \leq C_3 \qquad \forall \,\, B(x,2R_0) \subset \subset \R^{N}_{+},
	\end{equation}
where $C_3=C_3(N,f,R_0)>0$ is a constant.\\
From \eqref{bound-L-1-stab-out-10} and \eqref{estimate-any-p}-\eqref{estimate-f(u)-any-p>1} with $R=R_0$, we deduce that 
\begin{equation}\label{estimate-u-f(u)-any-p}
\begin{split}
   &\forall \, p>1 \quad \exists C_4=C_4(N,p,R_0,f,\theta)>0 \quad \text{such that} \\ 
   & \|u\|_{L^p({B(x,\frac{R_0}{4}}))} + \|f(u)\|_{L^p({B(x,\frac{R_0}{4}}))} \leq C_4 
   \qquad \forall \,\, B(x,2R_0) \subset \subset \R^{N}_{+}. \\
\end{split}   
\end{equation}

Then, by standard interior elliptic estimates, we have 
\begin{equation}\label{standard-estimate-u-f(u)-any-p}
\begin{split}
   &\forall \, p>1 \quad \exists C_5=C_5(N,p,R_0)>0 \quad \text{such that} \\ 
   & \|u\|_{W^{2,p}({B(x,\frac{R_0}{8}}))} \leq C_5 
   \Big [\|u\|_{L^p({B(x,\frac{R_0}{4}}))} + \|f(u)\|_{L^p({B(x,\frac{R_0}{4}}))} \Big]
   \qquad \forall \,\, B(x, 2R_0) \subset \subset \R^{N}_{+}. \\
\end{split}   
\end{equation}
Now we choose any $p>\frac{N}{2}$ in the latter and we apply \eqref{estimate-u-f(u)-any-p} and Sobolev embedding to get that 
\begin{equation}
    \|u\|_{L^{\infty}({B(x,\frac{R_0}{8}}))} \leq C_5 
    \qquad \forall \,\, B(x,2R_0) \subset \subset \R^{N}_{+},
\end{equation}
where $C_5=C_5(N,p,R_0, f, \theta)>0$ is a constant.\\
From the latter we deduce that $u$ is bounded on the half-space $\{x \in \R^{N}_{+} : x_N > 4R_0 \}$. The monotonicity of $u$ then implies that $u$ is bounded on $\R^{N}_{+}$. The desired conclusion the follows by invoking Theorem 4 in \cite{far-dup} (see Theorem H in the Introduction). 

The last statement of the theorem follows exactly as in the last part of the proof of Theorem \ref{symmetry-half-space-u-increase}. 

\end{proof}

\medskip

\begin{proof}[Proof of Corollary \ref{symmetry-half-space-u-bound-strip}] 
By the strong maximum principle, either $ u \equiv 0 $ (and we are done) or $u > 0 $ and $\frac{\partial u}{\partial x_{N}} > 0$ in $\R^{N}_{+}$. When $f(0) \geq 0$, any positive solution to \eqref{Poisson}, which is also bounded on finite strips, is strictly monotone in the $x_N$ direction (see Theorem 3.1 of \cite{far1}). The desired conclusions then follow by applying Theorem \ref{symmetry-half-space-u-increase} (resp. Theorem \ref{symmetry-half-space-u-increase2} or Theorem \ref{symmetry-half-space-u-increase3}) \end{proof}

\medskip

\begin{proof}[Proof of Theorem \ref{symmetry-half-space-superlinear}]
    For any $R\geq1$ and any $ x \in \R^{N}_{+}$ set $u_{R}(x):=u(Rx)$ and $B_R = B(0,R)$. Then, $u_{R}$ is a stable solution to  
\begin{equation}\label{Poisson-half-space-R}
	\left\{
	\begin{array}{ccl}
		-\Delta u_R= R^2f(u_R) & \text{in}&  \R^{N}_{+},\\
		u_R\geq 0 & \text{in} & \R^{N}_{+},\\
		u_R=0 & \text{on} & \partial \R^{N}_{+}.
	\end{array}
	\right.
\end{equation}     
We can therefore apply Proposition 5.5 in \cite{crfs} to \eqref {Poisson-half-space-R} to infer that
	\begin{equation}\label{estimates-gradient-superlinear}
		\disp\int_{B_{R/2}\cap \R^{N}_{+}}|\nabla u |^{2} \leq 
		C' R^{N-2}\Big( \frac{1}{R^{N}}\disp\int_{B_{R}\cap \R^{N}_{+}}u \Big)^{2},
	\end{equation}
where $C'>0$ is a dimensional constant. 

Now, since $f$ is non-negative and $R\geq1$, we observe that $u_{R}$ also solves 
\begin{equation*}
	-\Delta u_{R}=R^{2}f(u_{R})\geq f(u_{R}) \quad \text{in } \R^{N}_{+},
\end{equation*}
thus, in view of the superlinearity condition \eqref{superlinear1} on $f$, we can apply Theorem 3 in \cite{sirakov} to $u_{R}$ with $\Omega=B_{2}\cap \R^{N}_{+}$, $T=\partial \R^{N}_{+} \cap B_{2}$, $\omega=\Omega'=B_{1} \cap \R^{N}_{+}$, $A^{1}=Id$, $\lambda=1$, $q=2N$, $ \Lambda=\omega_{N}^{1/2N}$, $B=H=0$, $\varepsilon=1/2$, $p=1$, to get 
\begin{equation*}
	\disp\int_{B_{R}\cap \R^{N}_{+}} u \leq CR^{N}, \qquad \, \forall \,\,  R\geq1,
\end{equation*}
where $C=C(N,f)>0$ is a constant.

By plugging the latter into \eqref{estimates-gradient-superlinear} and recalling that $N \leq 4$, we obtain 
\begin{equation*}
      \disp\int_{B_{R}\cap \R^{N}_{+}} |\nabla u|^{2} \leq C'' R^{2}, \qquad \, \forall \,\,  R\geq1,
\end{equation*}
for some positive constant $C''$. Hence, according to Theorem \ref{two-dimensionnal-half-space}, up to a rotation in the first $N-1$ variables, we have that $u(x)=u_{0}(x_{N-1},x_{N}),$ where $u_0 \in C^{2}(\overline{\R^{2}_{+}})$ satisfies  
	\begin{equation*}
	\left\{
	\begin{array}{ccl}
		-\Delta u_{0}= f(u_{0}) & \text{in}&  \R^{2}_{+},\\
		u_{0} \geq 0 & \text{in} & \R^{2}_{+},\\
		u_{0}= 0 & \text{on} & \partial \R^{2}_{+}.
	\end{array}
	\right.
\end{equation*} 
Since $f \not\equiv 0$ by \eqref{superlinear1} and $f(0)\geq 0$, Theorem 2.5 in \cite{fs2} implies that $u_0 \equiv 0$, which in turn implies $u \equiv 0$ and $f(0) =0$. Furthermore, by plugging $u = 0$ into the stability condition we immediately get that $f'(0) \leq 0$, and the desired
conclusion $f'(0)=0$ then follows from the monotonicity assumption on $f$. 
\end{proof}

\bigskip

\begin{proof}[Proof of Theorem \ref{symmetry-half-space-N=3,4-dirichlet}]
By proceeding as in the proof of Theorem \ref{symmetry-half-space-superlinear} we immediately get
\begin{equation}\label{estimates-gradient-superlinear}
		\disp\int_{B_{R/2}\cap \R^{N}_{+}}|\nabla u |^{2} \leq 
		C' R^{N-2}\Big( \frac{1}{R^{N}}\disp\int_{B_{R}\cap \R^{N}_{+}}u \Big)^{2}, \qquad \, \forall \,\,  R\geq1,
	\end{equation}
where $C'>0$ is a dimensional constant. \\
Next we use the growth assumption \eqref{grand-O-half-space} to obtain that 
\begin{equation}
		\begin{split}
			 \Big( \frac{1}{R^{N}}\disp\int_{B_{R}\cap \R^{N}_{+}}u \Big)^{2}&\leq C(N) R^{2}\ln R, \quad \text{for } R >> 1, 
		\end{split}
	\end{equation}
where $C(N)>0$ is a dimensional constant. \\
The desired conclusion then follows by applying Theorem \ref{two-dimensionnal-half-space} of the present paper and Theorem 2.5 in \cite{fs2} exactly as in the last part of the proof of Theorem  \ref{symmetry-half-space-superlinear} (recall that $f \not\equiv 0$ by assumption).
\end{proof}

\bigskip

\begin{proof}[Proof of Theorem \ref{symmetry-half-space-f-convex-increasing}] 
By the strong maximum principle, either $u \equiv 0$ or $u>0$ in $\R^{N}_{+}$. The desired conclusion will follow if we rule out the latter case. For contradiction, assume that $u>0$ in $\R^{N}_{+}$. Since $f$ is a non-negative, convex function with $f(0)=0$, the function $t \mapsto \frac{f(t)}{t}$ is non-negative and non-decreasing on $(0, +\infty)$, so $\lim\limits_{t \to + \infty}\frac{f(t)}{t} \in [0,+\infty] $. Moreover, since $f \not \equiv 0$, we necessarily have $\lim\limits_{t \to + \infty}\frac{f(t)}{t} >0 $. 
Next we distinguish two cases:\\
\textit{Case 1 :} $\lim\limits_{s \to + \infty} \frac{f(s)}{s} = + \infty $.

By assumption, there exist a point $\bar x = ({\bar x}_1,...,{\bar x}_N) \in \R^{N}_{+}$ and an open ball $ B(\bar x, \bar R) \subset \subset \R^{N}_{+}$ such that $u$ is stable in $\R^{N}_{+} \setminus \overline{B(\bar x, \bar R)}$. If we set $\bar c =\frac{{\bar x}_N - \bar R} {4}>0$, then :


(i) for any $x \in \{x \in \R^{N}_{+} : x_N > \bar c \} \setminus \overline{B(\bar x,\bar R + \bar c)}$ we have that $B(x,\bar c) \subset \subset \R^{N}_{+}$ and that $B(x, \frac{{\bar c}} {2}) \subset \R^{N}_{+} \setminus \overline{B(\bar x, \bar R)}$ and therefore $u$ is stable in $B(x, \frac{{\bar c}} {2})$

and

(ii) $u$ is stable in the strip $\{x \in \R^{N}_{+} : x_N < 2 \bar c \}$.

When we are in the situation described in item (i), an application of Theorem 1.2 in \cite{crfs} provides 
\begin{equation}\label{estimate-alpha-holder-f-superlinear}
  	\|u\|_{L^{\infty}({B(x,\frac{{\bar c}} {4}}))} 
	\leq C_1 \|u\|_{L^{1}(B(x,\frac{{\bar c}} {2}))},
\end{equation}
where $C_1=C_1(N, \bar c)>0$ is a constant.\\
Since $f$ satisfies $\lim\limits_{s \to + \infty} \frac{f(s)}{s} = + \infty $ and $B(x,\bar c) \subset \subset \R^{N}_{+}$, we can apply item (ii) of Lemma \ref {lemma-bound-L-1} to get that
\begin{equation}\label{estimates-L1-f-superlinear}
	\disp\int_{B(x,\frac{{\bar c}} {2})} u \leq C_2,
\end{equation}
where $C_2 = C_2(f, N, \bar c)$ is a positive constant. \\
From \eqref{estimate-alpha-holder-f-superlinear} and \eqref{estimates-L1-f-superlinear} we obtain that $u \leq C_1 C_2$ on $ \{x \in \R^{N}_{+} : x_N \geq \bar c \} \setminus \overline{(B(\bar x,\bar R + \bar c)}$ and so, by the continuity of $u$ on $\overline{B(\bar x,\bar R + \bar c)}$, we deduce that $u$ is bounded in the closed half-space $ \{x \in \R^{N}_{+} : x_N \geq \bar c \}$.

On the other hand, for any $ x \in \{x \in \R^{N}_{+} : x_N < \bar c \}$, we have that $x \in B( (x',0), \bar c) \cap \R^{N}_{+}$. Since $u$ is stable in $\{x \in \R^{N}_{+} : x_N < 2 \bar c \}$ (see item (ii) above), we can apply Theorem 6.1 in \cite{crfs} to obtain that 
\begin{equation}\label{estimate-alpha-holder-boundary-f-superlinear222}
	\|u\|_{L^{\infty}(B((x',0), \bar c) \cap \R^{N}_{+})} \leq C_{3} \|u\|_{L^{1}(B((x',0),2\bar c)\cap \R^{N}_{+})},
	\qquad \, \forall \,\,  x' \in \partial \R^{N}_{+},
\end{equation}
where $C_3=C_3(N, \bar c)>0$ is a constant. \\
Since $f$ satisfies $\lim\limits_{s \to + \infty} \frac{f(s)}{s} = + \infty $, we can apply Theorem 3 in \cite{sirakov} to get that \begin{equation}\label{estimates-L1-f-superlinear-boundary-bbis}
	\disp\int_{B((x',0),2\bar c)\cap \R^{N}_{+}} u \leq C_4, \qquad \, \forall \,\,  x' \in \partial \R^{N}_{+},
\end{equation}
where $C_4 = C_4(f, N, \bar c)$ is a positive constant. \\
From \eqref{estimate-alpha-holder-boundary-f-superlinear222} and \eqref{estimates-L1-f-superlinear-boundary-bbis} we infer that $u$ is bounded on the strip $ \omega := \{x \in \R^{N}_{+} : x_N < \bar c \}$ and therefore $u$ is bounded on $\overline{\R^{N}_{+}}$.
Since $f(0)\geq 0$ and $u>0$ is bounded, Theorem 3.1 in \cite{far1} ensures that $u$ is strictly increasing in the $x_N$-direction. Therefore, $u$ is stable in $\R^{N}_{+}$ and we can apply once again Theorem 6.1 in \cite{crfs} to obtain that, for any $R>1$ and for any $x,y \in B(0,\frac{R}{2}) \cap \R^N_{+}$,
\begin{equation*}
	|u(x)-u(y)| \leq C_4 \|u\|_{L^{\infty}(\R^{N}_{+})} \frac{|x-y|^{\alpha}}{R^{\alpha}},
\end{equation*}
where $C_4>0$ is a dimensional constant. \\ 
By letting $R \to +\infty$, we obtain that $u(x)=u(y)$ for any $x,y \in \R^{N}_{+}$, which entails $u\equiv 0$ in $\R^{N}_{+}$, since $u=0$ on $\partial \R^{N}_{+}$. A contradiction with $u>0$ in $\R^N_{+}$.

\smallskip

\textit{Case 2 :} $\lim\limits_{t \to + \infty} \frac{f(t)}{t} \in (0, + \infty)$. 

In this case it is easily seen that $f$ is globally Lipschitz-continuous and Corollary 1.3 in \cite{BCN} implies that $u$ is strictly increasing in the $x_{N}$-direction. Hence, $u$ is stable on $\R^{N}_{+}$ and so Theorem 1.2 in \cite{crfs} tell us that
\begin{equation}
    \|u\|_{L^{\infty}({B(x,\frac{R_0}{2}}))} \leq C \|u\|_{L^{1}(B(x,R_0))} \qquad \text{for any} \quad B(x,R_0) \subset \subset \R^{N}_{+},
\end{equation}
where $C=C(N,R_0)>0$ is constant.\\
From the latter and item (i) of Lemma \ref{lemma-bound-L-1} we deduce that $u$ is bounded on the half-space $\{x \in \R^{N}_{+} : x_N > 2R_0 \}$, therefore, the monotonicity of $u$ then implies that $u $ is bounded on $\R^{N}_{+}$. We can then proceed as in the previous case to achieve a contradiction.   

So far, we have proved that $u\equiv 0$ on $\R^{N}_{+}$, the identities $f(0) = f'(0) = 0$ then follow as in the proof of Theorem \ref{symmetry-half-space-superlinear}. 
\end{proof}

\begin{proof}[Proof of Theorem \ref{symmetry-half-space-f-convex-increasing-croiss-contr}]
By the strong maximum principle, either $u \equiv 0$ or $u>0$ in $\R^{N}_{+}$. The desired conclusion will follow if we rule out the latter case. For contradiction, assume that $u>0$ in $\R^{N}_{+}$. \\ 
By assumption, there exist a point $\bar x = ({\bar x}_1,...,{\bar x}_N) \in \R^{N}_{+}$ and an open ball $ B(\bar x, \bar R) \subset \subset \R^{N}_{+}$ such that $u$ is stable in $\R^{N}_{+} \setminus \overline{B(\bar x, \bar R)}$. If we set 
$\bar c =\frac{{\bar x}_N - \bar R} {16}>0$, then :

(i) for any $x \in \{x \in \R^{N}_{+} : x_N > \bar c \} \setminus \overline{B(\bar x,\bar R + \bar c)}$ we have that $B(x,\bar c) \subset \subset \R^{N}_{+}$ and that $B(x, \frac{{\bar c}} {2}) \subset \R^{N}_{+} \setminus \overline{B(\bar x, \bar R)}$ and therefore $u$ is stable in $B(x, \frac{{\bar c}} {2})$

and

(ii) $u$ is stable in the strip $\{x \in \R^{N}_{+} : x_N < 8 \bar c \}$.

Now we observe that, without loss of generality, we can restrict ourselves to the 10-dimensional case. Indeed, if $u$ is stable in an 
open set $\omega \subset \R^N_+$, $N\leq 9$, then the function $\tilde{u}(y,x) := u(x)$ is stable in the open cylinder $\R^{10-N} \times \omega \subset \R^{10}_{+}$.\\
Hence, from now on, we set $ N=10 $. \\
When we are in the situation described in item (i), an application of Theorem 7.1 in \cite{crfs} provides 
\begin{equation}\label{estimate-any-p-N=10}
	\forall \, p>1 \quad \exists C_1=C_1(N,p,\bar c)>0 \quad \text{:} \quad \|u\|_{L^p({B(x,\frac{\bar c}{8}}))} \leq C_1 
    \|u\|_{L^{1}(B(x,\frac{\bar c}{2}))} \quad \forall \,\, B(x,\bar c) \subset \subset \R^{N}_{+}.
\end{equation}
Then, using the second condition in \eqref{superaff-u-croiss-contr} we get that 
\begin{equation}\label{estimate-f(u)-any-p}
\begin{split}
   &\forall \, p>1 \quad \exists C_2=C_2(N,p,\bar c,f,\theta)>0 \quad \text{such that} \\ &\|f(u)\|_{L^p({B(x,\frac{\bar c}{8}}))} 
   \leq C_2 \Big [1 + \|u\|^{\theta}_{L^{p \theta}(B(x,\frac{\bar c}{8}))} \Big]
   \qquad \forall \,\, B(x,\bar c) \subset \subset \R^{N}_{+}.
\end{split}   
\end{equation}
Since $f$ satisfies $\lim\limits_{s \to + \infty} \frac{f(s)}{s} = + \infty $ and $B(x,\bar c) \subset \subset \R^{N}_{+}$, we can apply item (ii) of Lemma \ref {lemma-bound-L-1} to get that
\begin{equation}\label{estimates-L1-f-superlinear-N=10}
	\disp\int_{B(x,\frac{{\bar c}} {2})} u \leq C_3,
\end{equation}
where $C_3 = C_3(f, N, \bar c)$ is a positive constant. \\
From \eqref{estimate-any-p-N=10}-\eqref{estimates-L1-f-superlinear-N=10} we deduce that 
\begin{equation}\label{estimate-u-f(u)-any-p-encore}
\begin{split}
   &\forall \, p>1 \quad \exists C_4=C_4(N,p,\bar c,f,\theta)>0 \quad \text{such that} \\ 
   & \|u\|_{L^p({B(x,\frac{\bar c}{8}}))} + \|f(u)\|_{L^p({B(x,\frac{\bar c}{8}}))} \leq C_4 
   \qquad \forall \,\, B(x,\bar c) \subset \subset \R^{N}_{+}. \\
\end{split}   
\end{equation}

Then, by standard interior elliptic estimates, we have 
\begin{equation}\label{standard-estimate-u-f(u)-any-p}
\begin{split}
   &\forall \, p>1 \quad \exists C_5=C_5(N,p,\bar c)>0 \quad \text{such that} \\ 
   & \|u\|_{W^{2,p}({B(x,\frac{\bar c}{16}}))} \leq C_5 
   \Big [\|u\|_{L^p({B(x,\frac{\bar c}{8}}))} + \|f(u)\|_{L^p({B(x,\frac{\bar c}{8}}))} \Big]
   \qquad \forall \,\, B(x,\bar c) \subset \subset \R^{N}_{+}. \\
\end{split}   
\end{equation}
Now we choose any $p > \frac{N}{2}$ in the latter and we apply \eqref{estimate-u-f(u)-any-p-encore} and Sobolev embedding to get that 
\begin{equation}
    \|u\|_{L^{\infty}({B(x,\frac{\bar c}{16}}))} \leq C_5 
    \qquad \forall \,\, B(x,\bar c) \subset \subset \R^{N}_{+},
\end{equation}
where $C_5=C_5(N,p,\bar c, f, \theta)>0$ is a constant.\\
From the latter we see that $u$ is bounded on $ \{x \in \R^{N}_{+} : x_N \geq \bar c \} \setminus \overline{(B(\bar x,\bar R + \bar c)}$ and so, by the continuity of $u$ on $\overline{B(\bar x,\bar R + \bar c)}$, we deduce that $u$ is bounded in the closed half-space $ \{x \in \R^{N}_{+} : x_N \geq \bar c \}$.

On the other hand, for any $ x \in \{x \in \R^{N}_{+} : x_N < \bar c \}$, we have that $x \in B( (x',0), \bar c) \cap \R^{N}_{+}$. 
Then we proceed as in the first part of the proof, but now we use boundary estimates. 
Since $u$ is stable in $\{x \in \R^{N}_{+} : x_N < 8 \bar c \}$ (see item (ii) above), we can apply Theorem 7.2 in \cite{crfs} to obtain that 
\begin{equation}\label{estimate-alpha-holder-boundary-f-superlinear22}
\begin{split}
    &\forall \, p>1 \quad \exists C_6=C_6(N,p,\bar c)>0 \quad \text{:} \quad \\
	& \|u\|_{L^{p}(B((x',0), 2\bar c) \cap \R^{N}_{+})} \leq C_6 \|u\|_{L^{1}(B((x',0),8\bar c)\cap \R^{N}_{+})}
	\qquad \, \forall \,\,  x' \in \partial \R^{N}_{+}
\end{split}	
\end{equation}
and using the second condition in \eqref{superaff-u-croiss-contr} as before, we infer that 
\begin{equation}\label{estimate-f(u)-any-p}
\begin{split}
   &\forall \, p>1 \quad \exists C_7=C_7(N,p,\bar c,f,\theta)>0 \quad \text{such that} \\ &\|f(u)\|_{L^{p}(B((x',0), 2\bar c) \cap \R^{N}_{+})}
   \leq C_7 \Big [1 + \|u\|^{\theta}_{L^{p}(B((x',0), 2\bar c) \cap \R^{N}_{+})} \Big]
   \qquad \, \forall \,\,  x' \in \partial \R^{N}_{+}. \\
\end{split}   
\end{equation}
Moreover, since $f$ satisfies $\lim\limits_{s \to + \infty} \frac{f(s)}{s} = + \infty $, we can apply Theorem 3 in \cite{sirakov} to get that \begin{equation}\label{estimates-L1-f-superlinear-boundary}
	\disp\int_{B((x',0),8 \bar c)\cap \R^{N}_{+}} u \leq C_7, \qquad \, \forall \,\,  x' \in \partial \R^{N}_{+},
\end{equation}
where $C_7 = C_7(f, N, \bar c)$ is a positive constant. \\
From \eqref{estimate-alpha-holder-boundary-f-superlinear22}-\eqref{estimates-L1-f-superlinear-boundary} we deduce that 
\begin{equation}\label{estimate-au-bord-u-f(u)-any-p-aubord}
\begin{split}
   &\forall \, p>1 \quad \exists C_8=C_8(N,p,\bar c,f,\theta)>0 \quad \text{such that} \\ 
   & \|u\|_{L^{p}(B((x',0), 2\bar c) \cap \R^{N}_{+})} + \|f(u)\|_{L^{p}(B((x',0), 2\bar c) \cap \R^{N}_{+})} \leq C_8 
   \qquad \, \forall \,\,  x' \in \partial \R^{N}_{+}. \\
\end{split}   
\end{equation}

Then, by standard elliptic boundary estimates on half-balls, we have 
\begin{equation}\label{standard-estimate-u-f(u)-any-p-aubord}
\begin{split}
   &\forall \, p>1 \quad \exists C_9=C_9(N,p,\bar c)>0 \quad \text{such that} \\ 
   & \|u\|_{W^{2,p}(B((x',0),\bar c) \cap \R^{N}_{+})} \leq C_9
   \Big [\|u\|_{L^p (B((x',0), 2\bar c) \cap \R^{N}_{+}))} + \|f(u)\|_{L^p(B((x',0), 2\bar c) \cap \R^{N}_{+}))} \Big]
   \qquad \, \forall \,\,  x' \in \partial \R^{N}_{+}. \\
\end{split}   
\end{equation}
Now we choose any $p>\frac{N}{2}$ in the latter and we apply \eqref{estimate-au-bord-u-f(u)-any-p-aubord} and Sobolev embedding to get that 
\begin{equation}
    \|u\|_{L^{\infty}(B((x',0),\bar c) \cap \R^{N}_{+})} \leq C_{10} 
    \qquad \, \forall \,\,  x' \in \partial \R^{N}_{+},
\end{equation}
where $C_{10}=C_{10}(N,p,\bar c, f, \theta)>0$ is a constant.	

From the latter we deduce that $u$ is bounded on the strip $\{x \in \R^{N}_{+} : x_N < \bar c \}$ and therefore $u$ is bounded on $\R^{N}_{+}$. Then, Theorem 4 in \cite{far-dup} (see Theorem H in the Introduction) implies that $u=u_{0}(x_N)$, where $u_0 : [0,+\infty) \to [0,+\infty)$ is a bounded, positive and strictly increasing function satisfying the ODE $ -u_0^{''} =f(u_0)$ on $[0,+\infty)$. Let us prove that this is impossible. Indeed, by the first integral of the ode we get that also $u'$ is bounded, so that $f(\sup u_{0})=0$ and 
$ 0 \leq f(u_0(t)) \leq f(\sup u_{0}) =0 $ for any $t \geq 0$. 
Therefore, $ -u_0^{''} = 0$, $u_0>0$ on $(0,+\infty)$ and $u_0(0)=0$, so that $u_{0}(t)=\alpha t $ with $\alpha >0$, contradicting the boundedness of $u_0$.


So far, we have proved that $u\equiv 0$ on $\R^{N}_{+}$, the identities $f(0) = f'(0) = 0$ then follow as in the proof of Theorem \ref{symmetry-half-space-superlinear}. 
\end{proof}

\section{ Proofs of Theorems \ref{sym-diff-ineq-1}-\ref{sym-diff-ineq-4}
}

\begin{proof}[Proof of Theorem \ref{sym-diff-ineq-1}]
By Lemma \ref{lemma-profile-infty} we know that $ \overline{u} \in L^1_{loc}(\R^{N-1})$ and so $v = - \overline{u}$ belongs to 
$L^1_{loc}(\R^{N-1})$ and satisfies 
\begin{equation}\label{ineq-Poisson-profile-infty-ineq}
\begin{split}
		-\Delta &v \leq -f(\overline{u}) \leq 0 \qquad \text{in} \quad  \mathcal{D}^{'}(\R^{N-1}), \\
		& v \leq 0 \qquad \text{on} \quad \R^{N-1}.
\end{split}		
\end{equation}
Since $N-1 =1,2$, then $\overline u$ must be constant (see for instance Theorem 3.2.18 in \cite{DiP-Val} or Theorem 3.2.24 in \cite{Horm}). Therefore, $f(\overline u)=0$ and so $\overline u =0$. The monotonicity of $u$, then implies $u \equiv 0$ on $\R^{N}_{+}$.   
\end{proof}

\medskip

\begin{proof}[Proof of Theorem \ref{sym-diff-ineq-2}]
For $ N=2,3$ the conclusion follows from Theorem \ref{sym-diff-ineq-1}. When $ N \geq 4$, let us distinguish the case $p=1$ from the case $p>1$. When $p=1$, by the strong maximum principle either $u \equiv 0$, and we are done, or $u>0$ on $\R^{N}_{+}.$ Let us prove that the latter case does not occur. Indeed, if $u>0$, then by multiplying the differential inequality $-\Delta u \geq u$ by $\frac{\varphi^2}{u}$, where $\varphi \in C^1_c(\R^{N}_{+})$, and integrating by parts we get $ \int_{\R^{N}_{+}} \vert \nabla \varphi  \vert^2 \geq \int_{\R^{N}_{+}} \varphi^2$. The latter and the variational characterization of $\lambda_1$ (the first eigenvalue of $-\Delta$ with zero Dirichlet boundary condition) tell us that $\lambda_1(B) \geq 1$ for any open ball $B \subset \subset \R^{N}_{+}$. This is clearly impossible, since $\lambda_1(B)$ approaches zero when $B$ is chosen arbitrary large.

When $p>1$, Lemma \ref{lemma-profile-infty} tells us that $ \overline{u} \in L^p_{loc}(\R^{N-1})$ and 
\begin{equation}\label{ineq-Poisson-profile-infty-ineq}
		-\Delta \overline{u} \geq \overline{u}^p  \qquad \text{in} \quad  \mathcal{D}^{'}(\R^{N-1}).
\end{equation}
Since $ 1 < p \leq \frac{N-1}{N-3} = \frac{(N-1)}{(N-1) -2}$, we can apply Theorem 2.1 in \cite{Miti-Poho} to obtain that $\overline u$ is a constant (necessarilly equal to zero). The monotonicity of $u$, then implies $u \equiv 0$ on $\R^{N}_{+}$.   
\end{proof}

\medskip

\begin{proof}[Proof of Theorem \ref{sym-diff-ineq-3}]
By the strong maximum principle, either $u \equiv 0$, and we are done, or $u>0$. Again, by applying the strong maximum principle to $\frac{\partial u}{\partial x_{N}}$ we see that, either $\frac{\partial u}{\partial x_{N}} =0$, and we are done, or 
$\frac{\partial u}{\partial x_{N}}>0$ in $\R^{N}_{+}$. The desired conclusion will follow if we rule out the latter case. 
For contradiction, assume that $u>0$ and $\frac{\partial u}{\partial x_{N}}>0$ in $\R^{N}_{+}$. 
The latter implies that $u$ is stable on $\R^{N}_{+}$ and so Theorem 1.2 in \cite{crfs} tell us that
\begin{equation}\label{estimate-alpha-holder-1-ineq}
    \|u\|_{L^{\infty}({B(x,\frac{R_0}{2}}))} \leq C \|u\|_{L^{1}(B(x,R_0))} \qquad \text{for any} \quad B(x,R_0) \subset \subset \R^{N}_{+},
\end{equation}
where $C=C(N,R_0)>0$ is constant.\\
The assumption \eqref{f-superlinear-N=2-ineq} and \eqref{estimate-alpha-holder-1-ineq} allow us to apply item (i) of Lemma \ref{lemma-bound-L-1} to deduce that $u$ is bounded on the half-space $\{x \in \R^{N}_{+} : x_N > 2R_0 \}$. The monotonicity of $u$ then implies that $u$ is bounded on $\R^{N}_{+}$. As a consequence, we see that $\overline u$ is a bounded classical solution to  
$-\Delta \overline{u} = f(\overline{u})$ in $\R^{N-1}$, $ N-1 \leq 8$. We can therefore invoke Theorem 1 (or Theorem 2) in \cite{far-dup} to prove that $\overline{u}$ is constant. Therefore, $f(\overline u)= 0$ and so $\overline u \equiv 0$. The monotonicity of $u$, then implies $u \equiv 0$ on $\R^{N}_{+}$. A contradiction.

If we further assume \eqref{superaff-u-croiss-contr2}, then we can proceed as in the proof of Theorem \ref{symmetry-half-space-u-increase3} to get, once again, that $u$ is bounded on some half-space 
$\{x \in \R^{N}_{+} : x_N > c >0\}$. Then, we conclude exactly as in the previous case (note that $ N-1 \leq 9$). 
\end{proof}

\medskip

\begin{proof}[Proof of Theorem \ref{sym-diff-ineq-4}]
Let us first prove that either $u \equiv 0$ or $u$ depends only on the first $N-1$ variables. To this aim, as in the proof of Theorem \ref{sym-diff-ineq-3}, it is enough to rule out the case : $u>0$ and $\frac{\partial u}{\partial x_{N}}>0$ in $\R^{N}_{+}$. For contradiction, assume $u>0$ and $\frac{\partial u}{\partial x_{N}}>0$ in $\R^{N}_{+}$, then $u$ is stable on $\R^{N}_{+}$ and so by Lemma \ref{lemma-profile-infty} we have $\overline u \in L^p_{loc}(\R^{N-1})$ and $\vert \overline u -u_n \vert^p \leq \vert \overline u \vert ^p$ a.e. on $\R^{N}$. Hence, by the dominated convergence theorem we infer that $u_n \longrightarrow \overline{u}$ in $L^p_{loc}(\R^{N})$ and $\overline u$ solves $-\Delta \overline{u} = \overline{u}^p $ in $ \mathcal{D}^{'}(\R^{N-1})$. 

Next we prove that $\overline{u}\in H^1_{loc}(\R^{N-1})\cap L^{p+1}_{loc}(\R^{N})$ and that $\overline{u}$ is a stable solution to $ -\Delta \overline{u} = f(\overline{u})$ in $\mathcal{D}^{'}(\R^{N-1}).$

To this end, we recall that $u$ is stable on $\R^{N}_+$ and so, from Proposition 4 in \cite{far} (applied with $\gamma=1$) we get
    \begin{equation}\label{estimate-Sobolev-1-ineq4}
	\int_{B(z,R)}|\nabla u |^{2} + |u|^{p+1} \leq C' R^{N - 2\frac{p+1}{p-1}} \qquad \text{for any } B(z,2R) \subset \subset \R^{N}_{+},
	\end{equation}
where $C'>0$ is a  constant depending only on $N$ and $p$.

Set $e_N =(0,...,0,1) \in \R^N$ and pick any open ball $B(z,r) \subset \R^{N}$. Then, there exists an integer 
$\bar n$, depending only on $z$ and $r$, such that 
$ B(z+n e_N, 2r) \subset \subset \R^{N}_{+}$ for any $ n \geq \bar n$.  Hence, by \eqref{estimate-Sobolev-1-ineq4} we deduce that 
$$
\int_{B(z,r)}|\nabla u_n |^{2} + |u_n|^{p+1} = \int_{B(z+n e_N, r)} |\nabla u |^{2} + |u|^{p+1}
\leq C' r^{N - 2\frac{p+1}{p-1}}  \qquad \text{for any } n \geq \bar n.
$$
From the latter, and by a standard diagonal argument, 
we can find a subsequence, still denoted by $(u_n)$, and $v \in H^1_{loc}(\R^{N})\cap L^{p+1}_{loc}(\R^{N})$  such that
\begin{equation}\label{cv-faible-u-n-p+1}
\begin{split}
&u_n \longrightarrow v \qquad \text{weakly in} \quad H^1_{loc} (\R^{N}), \\
& u_n \longrightarrow v \qquad \text{weakly in} \quad L^{p+1}_{loc}(\R^{N}). 
\end{split}
\end{equation}
Since we already know that $u_n \longrightarrow \overline{u}$ in $L^1_{loc}(\R^{N})$, from \eqref{cv-faible-u-n-p+1} we immediately get that $ \overline{u} = v \in H^1_{loc}(\R^{N}) \cap L^{p+1}_{loc}(\R^{N})$. Hence, $\overline{u} \in H^1_{loc}(\R^{N-1})\cap L^{p+1}_{loc}(\R^{N-1})$, since $\overline{u}$ is independent of the variable $x_N$. The stability of $ \overline{u}$ then follows exactly as in the last part of the proof of Theorem \ref{symmetry-half-space-u-increase2}. 

By summarizing, we have proved that $\overline{u}\in H^1_{loc}(\R^{N-1}) \cap L^{p+1}_{loc}(\R^{N})$ is a stable solution to $ -\Delta \overline{u} = \overline{u}^p$ in $\mathcal{D}^{'}(\R^{N-1})$, then Theorem 1.1 in \cite{Dav-Dup-Far} yields that $ \overline{u} \in C^2(\R^{N-1}).$ To conclude, we apply Theorem 1 in \cite{far} to get that $\overline{u} = 0$. The monotonicity of $u$, then implies $u \equiv 0$ on $\R^{N}_{+}$. A contradiction. 

So far, we proved that either $u \equiv 0$ or $u$ depends only on the first $N-1$ variables. To complete the proof we need only to prove the claim in item $(ii)$. To this end, we recall that all smooth positive solutions to the critical equation $ -\Delta v = v^{\frac{K+2}{K-2}}$ in $\R^{K}$, with $ K \geq 3$, have been classified in \cite{CGS} and are given by  
\begin{equation}\label{Talenti-Aubin-funct}
		v(x)= (a + b \vert x-x_0 \vert^2)^{-\frac{K-2}{2}} \qquad x \in \R^{K},
\end{equation}
for some $a,b>0$ with $1=abK(K-2)$ and $x_0 \in \R^K$. 
\end{proof}

\section{Auxiliary results}

In this section we state and prove some auxiliary results. They will be used to prove the main results of this article, but they are also interesting in their own right. Let us begin with two preliminaries results concerning non-negative solutions to the differential inequalities in the half-space. The first result is part of folklore, and we state it in a form useful for our purposes, while the second one deals with monotone solutions.

\begin{lem}\label{lemma-bound-L-1}
	Assume $N \geq 2$ and let $u \in C^2(\R^{N}_{+})$ be a solution to 
\begin{equation}\label{ineq-Poisson-gen-0}	
	\left\{
	\begin{array}{ccl}
		-\Delta u \geq f(u) & \text{in}&  \R^{N}_{+},\\
		u \geq 0 & \text{in} & \R^{N}_{+},
	\end{array}
	\right.
\end{equation} 
where $f \in C^0([0,+\infty))$.

(i) If $f$ satisfies 
	\begin{equation}\label{superlinear}
	f(t)\geq At-B \qquad \text{for any } t>0,
\end{equation}
where $A>0$ and $B \geq 0$, then, there exist $R_0 = R_0(A,N)>0$ and $C_0 = C_0(A,B,N)>0$ such that for any $R \geq R_0$ and any open ball $B(z,2R) \subset \subset \R^{N}_{+}$ we have 
    \begin{equation}\label{bound-L-1}
		\int_{B(z,R)}u \leq C_0 R^N. 
	\end{equation}

(ii) If $f$ satisfies 
\begin{equation}\label{superlinear-encore}
	\lim\limits_{t \to + \infty} \frac{f(t)}{t}=+\infty,
\end{equation}
then, for any $R>0$ and any open ball $B(z,2R) \subset \subset \R^{N}_{+}$ we have
\begin{equation}\label{bound-L-1-encore}
		\int_{B(z,R)}u \leq C_1, 
	\end{equation}
where $C_1=C_1(N,f,R)>0$ is a constant.  
\end{lem}

\begin{proof} (i) Denote by $\phi_1$ the positive first eigenfunction of $-\Delta$ on the open unit ball $B(0,1)$, such that $\max_{B(0,1)} \phi_1 = 1$ and by $\lambda_1>0$ the corresponding first eigenvalue. For any ball $ B(z,2R) \subset \subset \R^{N}_{+}$ set $\phi_{R,z}(x) := \phi_1(\frac{x-z}{2R})$, then 
 \begin{equation*}
	\left\{
	\begin{array}{cll}
	    -\Delta \phi_{R,z} = \frac{\lambda_1}{4R^2}\phi_{R,z}& \text{in} & B(z,2R),\\
		\phi_{R,z}>0& \text{in} & B(z,2R),\\
		\phi_{R,z}=0 & \text{on} & \partial B(z,2R).
	\end{array}
	\right.
\end{equation*}
From \eqref{ineq-Poisson-gen-0}	and \eqref{superlinear} we deduce that 
	\begin{equation}\label{est-affineA-B}
\int_{B(z,2R)} (Au-B)\phi_{R,z} \leq  - \int_{B(z,2R)} \Delta u \phi_{R,z} \leq  - \int_{B(z,2R)} u \Delta \phi_{R,z} = \frac{\lambda_1}{4R^2} \int_{B(z,2R)} u \phi_{R,z}. 
    \end{equation}
Set $ R_0= \sqrt{\frac{\lambda_1}{2A}}>0$, then for any $R \geq R_0$ we get
\begin{equation}
	\int_{B(z,2R)} Au \phi_{R,z} \leq \frac{\lambda_1}{4R_0^2} \int_{B(z,2R)} u \phi_{R,z} + B \int_{B(z,2R)} \phi_{R,z}  = \frac{A}{2} \int_{B(z,2R)} u \phi_{R,z} + B \int_{B(z,2R)}\phi_{R,z},
    \end{equation}
so that 
\begin{equation}
\frac{2B}{A} \int_{B(z,2R)} \phi_{R,z} \geq \int_{B(z,2R)} u \phi_{R,z} \geq \inf_{B(z,R)} \phi_{R,z} \int_{B(z,R)} u = \inf_{B(0,\frac{1}{2})} \phi_{1}\int_{B(z,R)} u.
\end{equation}
From the latter we get 
\begin{equation}
\int_{B(z,R)} u   \leq  (\inf_{B(0,\frac{1}{2})} \phi_{1})^{-1} \frac{2B}{A} \int_{B(z,2R)} \phi_{R,z}  \leq (\inf_{B(0,\frac{1}{2})} \phi_{1})^{-1} \frac{2B}{A} \int_{B(z,2R)} 1, 
\end{equation}
which implies \eqref{bound-L-1}. 

(ii) For any $R>0$ such that $B(z,2R) \subset \subset \R^{N}_{+}$, we set $A = \frac{\lambda_1}{4R^2}+1$, then by \eqref{superlinear-encore} there exists $B=B(f,A)>0$ such that $f(t)\geq At-B$ for any $t>0.$ We can therefore proceed as in the previous step to get the estimates \eqref{est-affineA-B} which, in turns, gives 
\begin{equation}
\int_{B(z,2R)} u \phi_{R,z} \leq  \int_{B(z,2R)} B\phi_{R,z} 
\end{equation}
and thus,
\begin{equation}
\int_{B(z,R)} u   \leq  (\inf_{B(0,\frac{1}{2})} \phi_{1})^{-1} B \int_{B(z,2R)} \phi_{R,z}  \leq (\inf_{B(0,\frac{1}{2})} \phi_{1})^{-1} B \int_{B(z,2R)} 1. 
\end{equation}
The latter implies the desired conclusion. 
\end{proof}

\medskip

\begin{lem}\label{lemma-profile-infty}
	Assume $N \geq 2.$ Let $u \in C^2(\R^{N}_{+})$ be a solution of 
\begin{equation}\label{ineq-Poisson-gen}	
	\left\{
	\begin{array}{ccl}
		-\Delta u \geq f(u) & \text{in}&  \R^{N}_{+},\\
		u \geq 0 & \text{in} & \R^{N}_{+},\\
		\frac{\partial u}{\partial x_{N}} \geq 0 & \text{in} & \R^{N}_{+},
	\end{array}
	\right.
\end{equation} 
where $f \in C^0([0,+\infty))$ satisfies 
	\begin{equation}\label{superlinear-bis}
	f(t)\geq At-B \qquad \text{for any } t>0,
\end{equation}
where $A>0$ and $B \geq 0$. Then,

$(i) \, $ $\overline{u}(x') := \lim_{x_N \to + \infty} u(x',x_N) \in [0,+\infty]$ for any $x'\in \R^{N-1}$ and $\overline{u} : \R^{N-1} \mapsto [0,+\infty]$ belongs to $L^1_{loc}(\R^{N-1})$.  

$(ii) \, $	The sequence of functions $(u_n)_{n \geq 1}$ defined by 
$$
u_n(x) := \tilde u (x',x_N + n), \qquad \text{for any } x \in \R^{N},
$$
converges to $\overline{u}$ in $L^1_{loc}(\R^{N})$.\footnote{ \, Here, as usual, we have denoted by $\tilde u$ the extension of $u$ with the value $0$ outside $\R^{N}_{+}$. \\
Also note that, the same convergence result holds true if we replace $u_n$ by 
$u_n(x) := \tilde u(x',x_N + a_n)$, where $(a_n)_{n \geq 1}$ is any sequence of positive real numbers converging monotonically to $+ \infty.$}

$(iii) \, $ $f(\overline{u}) \in L^1_{loc}(\R^{N-1}) $ and $\overline{u}$ solves 
    \begin{equation}\label{ineq-Poisson-profile-infty}
		-\Delta \overline{u} \geq f(\overline{u}) \qquad \text{in} \quad  \mathcal{D}^{'}(\R^{N-1}).
	\end{equation}
\end{lem}



\begin{proof} We observe that $\tilde u$ is measurable and so, $(u_n)_{n \geq 1}$ is a sequence of measurable functions satisfying
\begin{equation}\label{measur-monot-u-n}
		0 \leq u_n \leq u_{n+1} \leq \overline{u} \qquad \text{on } \R^N,  \qquad \text{for any } n \geq 1.
	\end{equation}

Set $e_N =(0,...,0,1) \in \R^N$, pick any $r \geq R_0$, where $R_0$ is given by item (i) of Lemma \ref{lemma-bound-L-1}
and consider any open ball $B(z,r) \subset \R^N$. There exists an integer $\bar n$, depending only on $z$ and $r$, such that 
$ B(z+n e_N, 2r) \subset \subset \R^{N}_{+}$ for any $ n \geq \bar n$. Hence, by Lemma \ref{lemma-bound-L-1} we deduce that 
$$
\int_{B(z+n e_N, r)} u \leq C_0(f,N) r^N \qquad \text{for any } n \geq \bar n,
$$
and
$$
\int_{B(z, r)} u_n =\int_{B(z+n e_N, r)} u \leq  C_0(f,N) r^N  \qquad \text{for any } n \geq \bar n.
$$
The latter and \eqref{measur-monot-u-n} then imply that $(u_n)_{n \geq 1}$ is a well-defined sequence in $L^1_{loc}(\R^N)$ satisfying 
$$
\int_{B(z, r)} u_n \leq C_0(f,N) r^N  \qquad \text{for any } n \geq 1.
$$
Then, the monotone convergence theorem implies that $u_n \longrightarrow \overline{u}$ in $L^1_{loc}(\R^{N})$. This proves items $(i)$ and $(ii)$.



Now, for any $n \geq 1$ and any $\varphi \in C^{\infty}_c(\R^{N}), \varphi \geq 0,$ we have $ 0 \leq (f(u_n) + B)\varphi$, hence 
\begin{equation}\label{Fatou1}
\int_{} (f(\overline{u}) +B)\varphi \leq \liminf_{n \to \infty} \int_{} (f(u_n) + B)\varphi = \liminf_{n \to \infty} \int_{} f(u_n) \varphi + \int_{} B\varphi
\end{equation}
by Fatou's Lemma. From the latter we deduce that $f(\overline{u})\in L^1_{loc}(\R^{N})$ (and so also in $L^1_{loc}(\R^{N-1})$, since it depends only on the first $N-1$ variables) and 
\begin{equation}\label{Fatou1-bis}
\int_{} f(\overline{u})\varphi \leq \liminf_{n \to \infty} \int_{} f(u_n) \varphi. \\
\end{equation}
As above, pick an open ball $B(z,r)$ such that $ supp\, \varphi \subset B(z,r)$, then we have $ -\Delta u_n \geq f(u_n)$ in $B(z,r)$, for any $n \geq \bar n.$ This property and \eqref{Fatou1-bis} 1ead to 
\begin{equation}\label{Fatou2}
\begin{split}
\int f(\overline{u}) \varphi & \leq \liminf_{n \to \infty} \int_{} f(u_n) \varphi \leq  \liminf_{n \to \infty} - \int \Delta u_n \varphi = \liminf_{n \to \infty} - \int u_n \Delta \varphi \\
& = \lim_{n \to \infty} - \int u_n \Delta \varphi = - \int \overline{u} \Delta \varphi. 
\end{split}
\end{equation}
The latter tell us that $ -\Delta \overline{u} \geq f(\overline{u})$ in $ \mathcal{D}^{'}(\R^{N})$ and the desired conclusion then follows by recalling that $\overline{u}$ only depends on the first $N-1$ variables. 
\end{proof}

The following theorem extends a result by N. Dancer valid for bounded stable solutions (see item (i) of Theorem 3 in \cite{dan1}). 
We remove the boundedness assumption and also make a slightly weaker hypotheses on the energy growth of the solutions over large balls.  

\begin{thm}\label{two-dimensionnal-half-space}
	Let $N\geq 2$ and $u \in C^{2}(\overline{\R^{N}_{+}})$ be a stable solution of
	\begin{equation*}
		\left\{
		\begin{array}{ccr}
			-\Delta u= f(u) & \text{in}&  \R^{N}_{+},\\
			u=0 & \text{in} & \partial \R^{N}_{+},
		\end{array}
		\right.
	\end{equation*} 
	where $f \in C^{1}(\R)$. Assume that
	\begin{equation}\label{grand-O}
		\disp\int_{B(0,R) \cap \R^{N}_{+}} |\nabla u|^{2}= O(R^{2} \ln R) \quad \text{as } R \to + \infty.
	\end{equation}
	Then, either \\
(i) $u$ is a function depending only on $x_{N}$,
	
or

(ii) (after a rotation in the first $N-1$ coordinates) $u$ is a function depending only on $(x_{N-1},x_{N})$ and strictly monotone in the  $x_{N-1}$ direction.
\end{thm}

The proof of Theorem \ref{two-dimensionnal-half-space} makes use of the following proposition and of the subsequent corollary. 

\begin{prop}\label{thm-general}
	Assume $N\geq 2$. Let $V\in C^{0}(\overline{\R^{N}_{+}})$, $\overline{\lambda}\geq 0$ and $v \in C^{1}(\overline{\R^{N}_{+}})$ be a weak solution to
\begin{equation}\label{eq-v}
	\left\{
	\begin{array}{ccr}
		-\Delta v= Vv+\overline{\lambda}v & \text{in}& \R^{N}_{+},\\
		v>0 & \text{in} & \R^N_+,\\
		v=0 & \text{on} & \partial \R^{N}_{+}.
	\end{array}
	\right.
\end{equation} 
Let $w \in C^{1}(\overline{\R^{N}_{+}})$ be a weak solution to
\begin{equation}\label{eq-w}
	\left\{
	\begin{array}{ccr}
		-\Delta w= Vw & \text{in}&  \R^N_+,\\
		w\leq 0 & \text{on} & \partial \R^{N}_{+},
	\end{array}
	\right.
\end{equation}
such that
\begin{equation}\label{grand-O-w}
	\disp\int_{B(0,R) \cap \R^N_+} (w^{+})^{2}  =O (R^{2} \ln R) \qquad \text{as } R\to +\infty.
\end{equation} 
Then, 
\begin{equation*}
	\frac{w^{+}}{v}\equiv constant \quad \text{in } \R^N_+,
\end{equation*}
where $w^{+}=\disp\max(w,0)$.
\end{prop}

\bigskip

\begin{cor}\label{corollary-thm-general}
	Let $v, w, V$ and $\overline{\lambda}$  as in the statement of Proposition \ref{thm-general}. Assume in addition that $w=0$ on $\partial \Omega$. Then
	\begin{equation*}
		\frac{w}{v}=constant \quad \text{in } \R^N_+.
	\end{equation*} 
\end{cor}

\bigskip

Let us now prove the results stated above.

\begin{proof}[Proof of Theorem \ref{two-dimensionnal-half-space}]

For any $R>0$ set $B_R = B(0,R)$. As in the proof of Theorem 3 in \cite{dan1}, we prove the existence of a positive function $v \in C^{2}(\overline{\R^{N}_{+}})$ such that 
	\begin{equation*}
		\left\{
		\begin{array}{ccl}
			-\Delta v=f'(u)v+\overline{\lambda}v & \text{in}&  \R^{N}_{+},\\
			v=0 & \text{on} & \partial \R^{N}_{+},
		\end{array}
		\right.
	\end{equation*}
	with $\overline{\lambda}= \lim\limits_{R \to + \infty} \lambda_{R} \geq 0$ where $(\lambda_{R},\phi_{R})$ is a pair eigenvalue-eigenfunction of $-\Delta \phi_{R}-f'(u)\phi_{R}=\lambda_{R}\phi_{R}$ in $B_{R} \cap \R^{N}_{+}$. Note that, since $f'$ and $u$ are continuous, then $f'(u)$ is bounded on compact sets of $\overline{\R^{N}_{+}}$ and Harnack inequality up to the boundary (Theorem 1.4 in \cite{BCNduke}) ensures that $(\phi_{R})_{R \geq R_{0}}$ is uniformly bounded on $B_{R_{0}}\cap \R^{N}_{+}$ for any $R_{0}>0$.  
	
Now, for $i \in \{1,\cdots, N-1\}$ we set $w = \frac{\partial u}{\partial x_i}$, then $w \in C^{1}(\overline{\R^{N}_{+}})$ is a weak solution to 
\begin{equation*}
	\left\{
	\begin{array}{ccl}
		-\Delta w=f'(u)w & \text{in}&  \R^{N}_{+},\\
		w=0 & \text{on} & \partial \R^{N}_{+},
	\end{array}
	\right.
\end{equation*}
and, since \eqref{grand-O} is in force, we get
\begin{equation*}
	\disp\int_{B_{R}\cap \R^{N}_{+}} w^{2} \leq \disp\int_{B_{R}\cap \R^{N}_{+}}|\nabla w|^{2}=O(R^{2}\ln{R}) \qquad \text{as } R \to + \infty.
\end{equation*}
Hence, according to Corollary \ref{corollary-thm-general} we have
	\begin{equation*}
		\frac{\partial u}{\partial x_i}=C_{i}v \quad \text{in } \R^{N}_{+}, \quad \text{ for any } i \in \{1,\cdots, N-1\},
	\end{equation*}
for some constants $C_i$.

	By an orthogonal rotation of axis in $\R^{N-1}$ we may and do suppose that $C_{i}=0$ for any $i\in \{1,N-2 \}$.
	
	Now, either $C_{N-1}=0$, and $u$ is a function of $x_{N}$ only, or $C_{N-1} \neq 0$ and, since $v >0$, then 
	$\frac{\partial u}{\partial {x_{N-1}}} = C_{N-1} v > 0 $ if $C_{N-1}>0$ ($<0$ if $C_{N-1}<0$). 
\end{proof}

\bigskip

\begin{proof}[Proof of Theorem \ref{thm-general}]
	Set $\sigma:=\frac{w}{v}\in C^{1}(\R^{N}_{+})$. Since $v$ and $w$ satisfy \eqref{eq-v} and \eqref{eq-w}, then $\sigma$ is a weak solution 
	to $div(v^{2}\nabla \sigma)= \overline{\lambda} v^{2}\sigma$ in $\R^{N}_{+}$. Then, by interior elliptic regularity theory, we have that $\sigma \in H^{2}_{\text{loc}}(\R^{N}_{+})$ and 
	\begin{equation}\label{eq-sigma}
		div(v^{2}\nabla \sigma)= \overline{\lambda} v^{2}\sigma \quad \text{a.e. in } \R^{N}_{+}.
	\end{equation}
	
For any $\Psi \in C^{0,1}_c(\R^{N})$, let $B$ an open ball in $\R^{N-1}$ and $T>0$ such that $\text{supp}(\Psi) \subset B \times  (-T,T):=\mathcal{C}$. 
For any integer $m > \frac{1}{T} $, multiply \eqref{eq-sigma} by $\sigma^{+} \Psi^{2}$ and integrate by parts over $\Omega_{m}\cap \mathcal{C}$, where $\Omega_{m}:=\{(x',x_N) \in \R^N : x_{N}> \frac{1}{m} \}$, to obtain
	\begin{equation}\label{int-by-parts_Omega_{m}}
		-\disp\int_{\Omega_{m}} v^{2} (\nabla \sigma . \nabla ( \sigma^{+} \Psi^{2})) - \underbrace{\disp\int_{\mathcal{C} \cap \partial \Omega_{m}} \sigma^{+} \Psi^{2} \Big(v \frac{\partial w}{\partial x_N}-w \frac{\partial v}{\partial x_N}\Big)}_{I_{m}}d\mathcal{H}^{N-1} =\disp\int_{\Omega_{m}}\overline{\lambda} (w^{+})^{2}\Psi^{2} \geq 0.
	\end{equation}

Next we prove that $\lim_{m \to \infty} I_m =0$. To this end we observe that, 
\begin{equation*}
	\begin{split}
	I_{m}&=\disp\int_{\mathcal{C} \cap \partial \Omega_{m} } \sigma^{+} \Psi^{2} v \frac{\partial w}{\partial x_N}d\mathcal{H}^{N-1}-\disp\int_{\mathcal{C} \cap \partial \Omega_{m} } \sigma^{+} \Psi^{2} w \frac{\partial v}{\partial x_N}d \mathcal{H}^{N-1}\\
	&=\underbrace{\disp\int_{\mathcal{C} \cap \partial \Omega_{m} } w^{+} \Psi^{2}  \frac{\partial w}{\partial x_N}d \mathcal{H}^{N-1}}_{I_{1,m}}-
	\underbrace{\disp\int_{\mathcal{C} \cap \partial \Omega_{m}} \sigma^{+} \Psi^{2} w^{+} \frac{\partial v}{\partial x_N}d \mathcal{H}^{N-1}}_{I_{2,m}}.
	\end{split}
\end{equation*}

Let us first deal with the term $I_{1,m}$. In this case we have 
\begin{equation*}
	I_{1,m} = \disp \int_{B} \bigg ( w^{+} \Psi^{2}  \frac{\partial w}{\partial x_N} \bigg )(x', 1/m) dx' : = \disp \int_{B} 
	g_m (x') dx'
\end{equation*}
with $\vert g_m \vert \leq \|w\|_{C^{1}(\overline{\R_+ \cap \mathcal{C}})}^{2}   \| \Psi \|_{C^{0}(\overline{\R^N_+ \cap \mathcal{C}})}^{2} \textbf{1}_B \in L^{1}(B)$ and $\lim\limits_{m \to + \infty}g_m(x')=0$  for any $x'\in B$. Thus, the dominated convergence theorem implies that 
\begin{equation}\label{limit_I1m}
	\lim\limits_{m \to + \infty} I_{1,m}=0.
\end{equation}

Now, we focus on the term $I_{2,m}$. Since $\sigma^{+}$ is not defined up to the boundary, we can not proceed directly as above. Nevertheless, by Hopf's Lemma we have $\frac{\partial v}{\partial x_N}>0$ on $\partial \R^N_+$ and so, by a standard compactness argument there exists $ \varepsilon >0$ such that   
\begin{equation}
		\frac{\partial v}{\partial x_{N}}(x) \geq \varepsilon \quad \text{ in } \overline{\mathcal{C} \cap \{(x',x_n) \in \R^N_+ : x_{N}< \varepsilon \}},
\end{equation}
which, in turn, implies that 
\begin{equation}\label{stima-v-bordo}
		v(x) \geq \varepsilon x_N \quad \text{ in } \overline{\mathcal{C} \cap \{(x',x_n) \in \R^N_+ : x_{N}< \varepsilon \}}.
\end{equation}
Also, since $w^+ = 0$ on $\partial \R^N_+$ we have 
\begin{equation}\label{stima-v-bordo-bis}
		w^+(x) \leq \|w\|_{C^{1}(\overline{\R_+ \cap \mathcal{C}})} x_N \quad \text{ in } \overline{\mathcal{C} \cap \{(x',x_n) \in \R^N_+ : x_{N}< \varepsilon \}}.
\end{equation}
Now we estimate the term $I_{2,m}$ for $ m > \max \{\frac{1}{\varepsilon}, \frac{1}{T}\}$. Since 
\begin{equation*}
	I_{2,m} = \disp \int_{B} \bigg ( \sigma^{+} \Psi^{2} w^{+} \frac{\partial v}{\partial x_N} \bigg )(x', 1/m) dx' : = \disp \int_{B} h_m (x') dx'
\end{equation*}
and 
\begin{equation*}
\sigma^{+} (x', 1/m)= \frac{w^{+} (x', 1/m)}{v(x', 1/m)} \leq  \frac{\|w\|_{C^{1}(\overline{\R_+ \cap \mathcal{C}})} }{m} 
\frac{m}{\varepsilon} = \frac{\|w\|_{C^{1}(\overline{\R_+ \cap \mathcal{C}})} }{\varepsilon} \qquad \text{for any } x' \in B, 
\end{equation*}
where in the latter inequality we have used \eqref{stima-v-bordo}-\eqref{stima-v-bordo-bis}, we deduce that $h_m$ is bounded on $B$, independently on $m > \max \{\frac{1}{\varepsilon}, \frac{1}{T}\}$. Furthermore, $\lim\limits_{m \to + \infty}h_m(x')=0$  for any $x'\in B$, therefore $\lim\limits_{m \to + \infty} I_{2,m}=0$ by the dominated convergence theorem. Hence, $\lim\limits_{m \to + \infty}I_{m}=0$.

We also observe that $v^{2} (\nabla \sigma . \nabla ( \sigma^{+} \Psi^{2}))$ is bounded on $\R^N_+$ and zero a.e. outside $B \times (0,T)$, hence $ \int_{\Omega_{m}} v^{2} (\nabla \sigma . \nabla ( \sigma^{+} \Psi^{2})) \longrightarrow \int_{\R^N_+} v^{2} (\nabla \sigma . \nabla ( \sigma^{+} \Psi^{2}))$, as $ m \to + \infty$. 

Therefore, by letting $ m \to + \infty$ in \eqref{int-by-parts_Omega_{m}}, we get that
\begin{equation}\label{inequation-sigma}
-\disp\int_{\R^N_+} v^{2} (\nabla \sigma . \nabla ( \sigma^{+} \Psi^{2})) \geq 0.
\end{equation}	
Then, by Young inequality we obtain
\begin{equation*}
	\begin{split}
	\disp\int_{\R^{N}_{+}} v^{2} |\nabla \sigma^{+}|^{2} \Psi^{2} \leq -2 \disp\int_{\R^N_+} v^{2} \sigma^{+} \Psi (\nabla \sigma^{+}. \nabla \Psi) &\leq \disp\int_{\R^N_+} (v \Psi |\nabla \sigma^{+}|) (2v \sigma^{+} |\nabla \Psi |)\\
	& \leq  \frac{1}{2} \disp\int_{\R^N_+} v^{2}|\nabla \sigma^{+}|^{2}  \Psi^{2}  + 2\disp\int_{\R^N_+} v^{2} (\sigma^{+})^{2} |\nabla \Psi|^{2}, 
	\end{split}
\end{equation*}
which entails
\begin{equation}\label{inequalities-young}
	\disp\int_{\R^N_+} v^{2} |\nabla \sigma^{+}|^{2} \Psi^{2} \leq 4\disp\int_{\R^N_+} v^{2} (\sigma^{+})^{2} |\nabla \Psi|^{2}. 
\end{equation}
For $R>1$, let 
$\Psi_{R}$ be defined by 
\begin{equation*}
	\Psi_{R}(x):=\textbf{1}_{B_{\sqrt{R}}}(x)+\dfrac{2\ln(R/|x|)}{\ln R}\textbf{1}_{B_{R}\backslash B_{\sqrt{R}}}(x),
\end{equation*}
then
\begin{equation}\label{estimate-v2-nabla-sigma2}
	\disp\int_{ B_{R} \cap \R^N_+} v^{2}  (\sigma^{+})^{2} |\nabla \Psi_{R}|^{2} \leq \frac{4}{\ln{R}^{2}}\disp\int_{ (B_{R} \backslash B_{\sqrt{R}}) \cap \R^N_+} \frac{(w^{+})^{2}}{|x|^{2}}.
\end{equation}
Now we use the identity $\frac{1}{|x|^{2}}=\int_{|x|}^{R} \frac{2}{t^{3}} +\frac{1}{R^{2}}$ and the assumption \eqref{grand-O-w} to get
\begin{equation}\label{inequality-1-divided-by-x2}
	\begin{split}
		\disp\int_{(B_{R}\backslash B_{\sqrt{R}})\cap \R^{N}_{+}}\frac{(w^{+})^{2}}{|x|^{2}}&=\disp\int_{\sqrt{R}}^{R}\disp\int_{(B_{t}\backslash B_{\sqrt{R}})\cap \R^N_+} \frac{2(w^{+})^{2}}{t^{3}}dxdt +\disp\int_{(B_{R}\backslash B_{\sqrt{R}}) \cap \R^N_+} \frac{(w^{+})^{2}}{R^{2}}\\
		&\leq 2C \ln^{2}R +C \ln R.
	\end{split}
\end{equation}
where $C>0$ is a constant independent of $R$ (such that $\int_{ B_{R} \cap \R^N_+ } (w^{+})^{2} \leq C R^{2}\ln R$).\\
From the latter, \eqref{inequalities-young} and \eqref{estimate-v2-nabla-sigma2} we deduce that
\begin{equation}\label{doppia-disug-finale}
\disp\int_{\R^N_+} v^{2} |\nabla \sigma^{+}|^{2} \Psi_R^{2} \leq 4\disp\int_{\R^N_+} v^{2} (\sigma^{+})^{2} |\nabla \Psi_R|^{2} \leq 48 C \qquad \text{for } R>>1, 
\end{equation}
In particular, by letting $R \to \infty$, we get $\int_{\R^N_+} v^{2} |\nabla \sigma^{+}|^{2} \leq 48 C$, that is $v^{2} |\nabla \sigma^{+}|^{2} \in L^{1}(\R^N_+).$ \\
Now, from \eqref{inequation-sigma}, \eqref{doppia-disug-finale} and Cauchy-Schwarz inequality we obtain 
\begin{equation*}
	\begin{split}
		&\disp\int_{B_{R}\cap \R^N_+} v^{2} |\nabla \sigma^{+}|^{2} \Psi_{R}^{2} \leq -2 \disp\int_{\R^N_+} v^{2} \sigma^{+} \Psi_R (\nabla \sigma^{+}. \nabla \Psi_R) \\
		 &\leq 2 \Big( \disp\int_{(B_{R}\backslash B_{\sqrt{R}}) \cap \R^N_+}   v^{2} |\nabla \sigma^{+}|^{2} \Psi_{R}^{2} \Big)^{1/2} \Big( \disp\int_{(B_{R}\backslash B_{\sqrt{R}}) \cap \R^N_+} v^{2} (\sigma^{+})^{2} |\nabla \Psi_{R}|^{2}  \Big)^{1/2}\\
		 &\leq 2 \Big( \disp\int_{(B_{R}\backslash B_{\sqrt{R}}) \cap \R^N_+}   v^{2} |\nabla \sigma^{+}|^{2} \Psi_{R}^{2} \Big)^{1/2} \Big( \disp\int_{B_{R} \cap \R^N_+} v^{2} (\sigma^{+})^{2} |\nabla \Psi_{R}|^{2}  \Big)^{1/2}\\
		 & \leq  8\sqrt{C} \Big( \disp\int_{(B_{R}\backslash B_{\sqrt{R}}) \cap \R^N_+}   v^{2} |\nabla \sigma^{+}|^{2} \Psi_{R}^{2} \Big)^{1/2} \leq 8 \sqrt{C}\Big( \disp\int_{(B_{R}\backslash B_{\sqrt{R}}) \cap \R^N_+}   v^{2} |\nabla \sigma^{+}|^{2}  \Big)^{1/2}
	 \end{split}
\end{equation*}
Finally, since $v^{2}|\nabla \sigma^{+}|^{2} \in L^{1}(\R^N_+)$, we deduce from the dominated convergence theorem that $\int_{(B_{R}\backslash B_{\sqrt{R}}) \cap \R^N_+}   v^{2} |\nabla \sigma^{+}|^{2} \to 0$ as $R \to + \infty$, and so $\int_{\R^N_+}v^{2} |\nabla \sigma^{+}|^{2}=0$. Since 
 $v>0$ in $\R^N_+$, we infer that $\nabla \sigma^{+}=0$ almost everywhere in $\R^N_+$ and so $\frac{w^{+}}{v} = \sigma^{+} = constant$ in $\R^N_+$.
\end{proof}

\bigskip

\subsection{Notations}\label{Notations} \quad \\

\noindent  $\mathbb{R}^N_+ =\{x=(x', x_N) \in \R^{N-1} \times \R  \ | \ x_N>0\}$, the open upper half-space of $\R^N$.  

\bigskip

\noindent $ \vert \cdot \vert$ : the Euclidean norm. 

\bigskip

\noindent $B(x,R)$  : the Euclidean $N$-dimensional open ball of center $x$ and radius $R>0$. 

\bigskip

\noindent $B'(x',R)$  : the Euclidean $N-1$-dimensional open ball of center $x'$ and radius $R>0$. 

\bigskip

\noindent $B_R : = B(0 , R)$  and $B'_R := B'(0', R)$, where $0 =(0',0) \in \R^{N-1} \times \R$ is the origin of $ \R^N$.  

\bigskip

\noindent $\mathcal{H}^{N-1}$ : the $N-1$-dimensional Hausdorff measure. 

\bigskip

%

\noindent $Lip(X)$ : the set of globally Lipschitz-continuous functions on $X$. 

\bigskip

\noindent $Lip_{loc}(X)$ : the set of locally Lipschitz-continuous functions on $X$. 

\bigskip

\noindent $C^k(\overline{U})$ : the set of functions in $C^k(U)$ all of whose derivatives of order $ \leq k$ have continuous (not necessarily bounded) extensions to the closure of the open set $U$. 

\bigskip

\noindent $C^{0,\alpha}(\overline{U})$ : the vector space of bounded and globally $\alpha$-Hölder-continuous functions $h$ on the open set $U$ endowed with the norm :

\medskip

\qquad $\Vert h \Vert_{C^{0,\alpha}(\overline{U})} := \|h\|_{L^{\infty}(U)} + 
\left[ h \right]_{C^{0,\alpha}  (U)} :=  \sup_{x \in U} |h(x)| + \sup_{x,y \in U, x\neq y}
\frac{|h(x)-h(y)|}{|x-y|^{\alpha}} $ . 

\bigskip

\noindent $C^{k,\alpha}(\overline{U})$ : the vector space of functions in $C^k(U)$ all of whose derivatives of order $ \leq k$ belong to $C^{0,\alpha}(\overline{U})$, endowed with the norm :

\medskip

\hskip4truecm  $\Vert h \Vert_{C^{k,\alpha}(\overline{U})} := \sum_{0 \leq \vert \beta \vert \leq k} \Vert \partial^{\beta} h 
\Vert_{C^{0,\alpha}(\overline{U})}$ . 

\bigskip

\noindent $\mathcal{D}'(U)$ : the space of distributions on the open set $U$. 


%

\bibliography{plain.bst}

\end{document}